\documentclass[final]{siamltex}

\usepackage{cite}
\usepackage{comment}
\usepackage{amsmath}
\usepackage{graphicx}
\usepackage{fancyhdr}
\usepackage{pslatex}
\usepackage{epsfig,psfrag}
\usepackage{amsmath}
\usepackage{amssymb}
\usepackage{latexsym,pifont,color,comment}
\usepackage{stmaryrd}
\usepackage{esint}
\usepackage{algorithm,algorithmic}

\def\B{{\bf B}}
\def\A{{\bf A}}
\def\u{{\bf u}}
\def\En{{\mathcal E}}

\def\E{{\bf E}}
\def\n{{\bf n}}
\def\nv{{\bf n}}
\def\xv{{\bf x}}
\def\NF{{\mathcal F}}
\def\WS{{\mathcal W}}
\def\meq{{M_\text{eq}}}

\def\Bt{\widetilde{B}}

\def\deg{{\text q}}

\def\ent{U}
\def\dualent{U^{\star}}
\def\eflux{{\bf G}}
\def\dualeflux{{\bf G}^{\star}}

\def\beq{\begin{equation}}
\def\eeq{\end{equation}}

\def\bga{\begin{gather}}
\def\ega{\end{gather}}

\def\bal{\begin{align}}
\def\eal{\end{align}}

\def\Tm{{\mathcal T}}

\def\bmat{\begin{bmatrix}}
\def\emat{\end{bmatrix}}
\def\reals{{{\rm l} \kern -.15em {\rm R} }}

\newcommand{\al}[1]{\ensuremath{b_{i-\half \, j}^{1 \, (#1)}}}
\newcommand{\ar}[1]{\ensuremath{b_{i+\half \, j}^{1 \, (#1)}}}
\newcommand{\bl}[1]{\ensuremath{b_{i \, j-\half}^{2 \, (#1)}}}
\newcommand{\br}[1]{\ensuremath{b_{i \, j+\half}^{2 \, (#1)}}}

\newcommand{\aal}[1]{\ensuremath{b_{i-\half \, j \, k}^{1 \, (#1)}}}
\newcommand{\aar}[1]{\ensuremath{b_{i+\half \, j \, k}^{1 \, (#1)}}}
\newcommand{\bbl}[1]{\ensuremath{b_{i \, j-\half \, k}^{2 \, (#1)}}}
\newcommand{\bbr}[1]{\ensuremath{b_{i \, j+\half \, k}^{2 \, (#1)}}}
\newcommand{\ccl}[1]{\ensuremath{b_{i \, j \, k-\half}^{3 \, (#1)}}}
\newcommand{\ccr}[1]{\ensuremath{b_{i \, j \, k+\half}^{3 \, (#1)}}}

\def\half{\frac{1}{2}}
\def\halfsqt{\frac{1}{2\sqrt{3}}}
\def\halfsqft{\frac{1}{2\sqrt{15}}}
\def\halfsqtf{\frac{1}{2\sqrt{35}}}
\def\dxody{\frac{\Delta x}{\Delta y}}
\def\dyodx{\frac{\Delta y}{\Delta x}}
\def\dxodz{\frac{\Delta x}{\Delta z}}
\def\dzodx{\frac{\Delta z}{\Delta x}}
\def\dzody{\frac{\Delta z}{\Delta y}}
\def\dyodz{\frac{\Delta y}{\Delta z}}

\newtheorem{thm}{Theorem}
\newtheorem{remark}[theorem]{Remark}

\title{\bf High-Order Discontinuous Galerkin Finite Element Methods with Globally Divergence-Free 
Constrained Transport for Ideal MHD}

\author{James A. Rossmanith\thanks{Iowa State University, Department of Mathematics,
396 Carver Hall, Ames, IA 50011, USA ({\tt rossmani@iastate.edu})}}

\begin{document}

\maketitle

\begin{abstract}
The modification of the celebrated Yee scheme from the vacuum Maxwell equations
to magnetohydrodynamics (MHD) is often referred to as the constrained transport (CT)
approach. Constrained transport can be viewed as a sort of predictor-corrector
method for updating the magnetic field, where a magnetic field value is first predicted
by a method that does not exactly preserve the divergence-free condition on the
magnetic field, followed by a correction step that aims to control these divergence
errors. This strategy has been successfully used in conjunction
with a variety of shock-capturing methods including WENO (weighted
essentially non-oscillatory), central, and wave propagation schemes.
In this work we show how to extend the basic CT framework in the context of
the discontinuous Galerkin (DG) finite element method on both
2D and 3D Cartesian grids. We first
 review the entropy-stability theory for semi-discrete DG discretizations
 of ideal MHD, which rigorously establishes the need for
 a magnetic field that satisfies the following conditions:
 (1) the divergence of the magnetic field is zero on each element, and (2)
 the normal components of the magnetic field are continuous across all element edges (faces in 3D).
 In order to achieve such a globally divergence-free magnetic field,  we introduce a novel constrained transport
 scheme that is based on two main ingredients: (1) we introduce an element-centered magnetic vector potential
 that is updated via a discontinuous Galerkin scheme on the induction equation; and (2) we define
 a mapping that takes element-centered magnetic field values (i.e., the predicted magnetic field)
 and element-centered magnetic vector potential values and creates on each edge  (face in 3D)
 a high-order representation of the normal component of the magnetic field; this representation is then
 mapped back to the elements to create a globally divergence-free element-centered representation
 of the magnetic field.   For problems with shock waves, we make use of so-called moment-based limiters
 to control oscillations in the conserved quantities. The resulting method is applied to standard test cases for ideal MHD.
\end{abstract}



\section{Introduction}
{Plasma} is often referred to as the fourth state of matter after solid, liquid, and gas,
and consists of a mixture of interacting charged particles. 
Macroscopic features of a quasi-neutral plasma can often be accurately
modeled through {magnetohydrodynamic} (MHD) models that
track only macroscopic quantities such as the total mass density, center-of-mass
momentum density, and total energy density
(see standard plasma physics textbooks such as Chapter 4 of Gombosi \cite{book:Go98}).
The {ideal} magnetohydrodynamic equations further assume that the flow is
 inviscid and that the plasma is a perfect conductor (i.e., zero resistivity). 

The ideal MHD system can be written as a system of hyperbolic
conservation laws, where the conserved quantities 
are mass density, momentum density, energy density, and the magnetic field.
Furthermore, this system is equipped with an entropy inequality that
features a convex scalar entropy and a corresponding
entropy flux. Indeed, the scalar entropy, with  some
help from the fact that the magnetic field is divergence-free, 
can be used to define entropy variables in which the 
MHD system is in symmetric hyperbolic form \cite{article:Barth05,article:Go72}.

As has been noted many times in the literature (e.g., 
Brackbill and Barnes \cite{article:BrBa80}, Evans and Hawley \cite{article:EvHa88},
and T\'oth \cite{article:To00}),
numerical methods for ideal MHD must in general satisfy (or at least control)
some discrete version of the divergence-free condition on the magnetic field:
\begin{equation}
\label{eqn:divBeq0}
\nabla \cdot \B = 0.
\end{equation}
Failure to accomplish this generically leads to a nonlinear numerical instability,
which often leads to negative pressures and/or densities.
Starting with the paper of Brackbill and Barnes \cite{article:BrBa80} in 1980,
several approaches for controlling errors in $\nabla \cdot \B$ have
been proposed. An in-depth review of many of these methods can
be found in T\'oth \cite{article:To00}. 

The constrained transport (CT) approach for ideal MHD was introduced by Evans and Hawley
\cite{article:EvHa88}. The method is a modification of Yee's method \cite{article:Ye66} for electromagnetic 
wave propagation, and, at least in its original formulation, introduced staggered magnetic and electric fields.
This approach can also be viewed as a kind of
predictor-corrector scheme for the magnetic field. Roughly speaking, the idea is
to compute all of the conserved quantities with a ``standard'' finite difference, finite volume, or finite element method. 
This step produces the {\it predicted} magnetic field values.
From these computed quantities one
then constructs an approximation to the electric field through the ideal Ohm's law.
This electric field can then be used to update the magnetic vector potential,
which in turn, can be used to compute a divergence-free magnetic field. This
step produces the {\it corrected} magnetic field values.

The main advantages of this approach compared to other approaches 
for solving the ideal MHD equations are that (1) there is no elliptic solve
such as in projection methods (e.g., see T\'oth \cite{article:To00} and Balsara and Kim \cite{article:BaKi04}), 
and (2) there are no free parameters to choose such as in the hyperbolic
divergence-cleaning technique \cite{article:Ded01a}.
The main disadvantages of this approach are that (1) additional variables must be 
stored and updated  (i.e., the magnetic potential), and (2) a staggered description
typical of mixed finite element methods is often required\footnote{There are CT methods that
avoid a staggered magnetic field (see for example \cite{article:To00,article:Ro04b,article:HeRoTa11,
article:HeRoTa13,article:ChRoTa14,article:Ro04a}).}.

Since the introduction of the original CT framework by Evans and Hawley \cite{article:EvHa88}, several variants and modifications 
have been introduced, including the work of Balsara \cite{article:Ba04},  Balsara and Spicer \cite{article:BaSp99a}, 
Christlieb et al. \cite{article:ChRoTa14}, Dai and Woodward \cite{article:DaWo98}, 
Fey and Torrilhon \cite{article:FeTo03}, 
Helzel et al. \cite{article:HeRoTa11,article:HeRoTa13},
Londrillo and Del Zanna \cite{article:LoZa04},
Rossmanith \cite{article:Ro04b}, Ryu et al. \cite{article:Ryu98}, and De Sterck \cite{article:De01b}. An overview of many of these approaches, as well as the
introduction of a few more variants, can be found in
T\'oth \cite{article:To00}. 

In this work we show how to extend the basic constrained transport
 framework in the context of the discontinuous Galerkin finite element method (DG-FEM) 
 on both 2D and 3D Cartesian meshes. 
The method advocated in this work makes use of two key ingredients:
(1) the introduction of a magnetic vector potential, which in our case will be represented
in the same finite element space as the conserved variables, and (2)
the use of a particular divergence-free
reconstruction of the magnetic field, which makes use of the magnetic vector potential
and the predicted magnetic field.
The divergence-free reconstruction advocated in this work is slight modification of the reconstruction
method that has been used in other work on ideal MHD. Li, Xu, and Yakovlev \cite{article:LiXuYa11} made
use of this reconstruction in the context of a 2D central DG scheme. Balsara
\cite{article:Ba04} made use this reconstruction in the context of finite volume schemes
and adaptive mesh refinement. In fact, the divergence-free reconstruction of the magnetic
field is directly related to ideas from discrete exterior calculus and, in particular,
to Whitney forms (e.g., see Bossavit \cite{article:Bossavit88}). Despite these connections, we do need
the full machinery of discrete exterior calculus and Whitney forms in order to
develop the proposed numerical scheme.
The novel aspect of this work is that we make direct use of a magnetic vector potential,
thus following in the footsteps of the methods developed in
\cite{article:ChRoTa14,article:HeRoTa11,article:HeRoTa13,article:LoZa04,article:Ro04b,article:De01b}.
An advantage of our approach is that extension from 2D to 3D is straightforward.
 
The paper begins with a description of the ideal MHD equations in \S \ref{sec:equations}.
In  \S \ref{sec:dg_method} we describe the basic DG method, both on 2D and 3D Cartesian meshes.
In \S \ref{sec:dg_stability} we review the entropy stability theory for semi-discrete DG-FEM,
including a discussion of the relevant theorem of Barth \cite{article:Barth05} that rigorously establishes the need for
divergence-free magnetic fields in the discretization of MHD.
In \S \ref{sec:ct_scheme} we introduce a novel constrained transport scheme that is based on achieving a globally divergence-free magnetic field via the use
of a combination of techniques from mixed finite element methods, as well as through the use of a magnetic potential. The magnetic potential is updated through an appropriate high-order discretization of the induction equation, and a globally divergence-free magnetic field is then constructed through a discrete curl operation on the magnetic potential. In particular, the discrete curl of the magnetic potential is defined through a high-order construction of the normal components of the magnetic field on each element edge, followed by a reconstruction step to obtain a globally-defined element-centered definition of the magnetic field.  This newly proposed method is closely related to the
central DG scheme for ideal MHD developed by Li, Xu, and Yakovlev \cite{article:LiXuYa11}, although the new
scheme does not require the use of central DG, nor does it require storing the solution on both a primal
and a dual grid. For these reasons, extension to 3D is straightforward.
The resulting scheme is applied to several numerical test cases in \S \ref{sec:numex}.
All of the methods presented in this work have been implemented in the {\sc DoGPack} software package \cite{dogpack}
and will be made freely available on the web. All the visualization in this work
has been done with the freely available package {\sc Matplotlib} \cite{article:matplotlib}.

 \section{Ideal magnetohydrodynamic equations}
\label{sec:equations}
The ideal magnetohydrodynamic (MHD) equations are a classical model
from plasma physics that describe the macroscopic evolution of a 
quasi-neutral two-fluid plasma system. Under the quasi-neutral
assumption, the two-fluid equations can be collapsed
into a single set of fluid equations for the total mass density, center-of-mass momentum density, and 
total energy density of the system. The resulting equations can be written in the following 
form\footnote{We use boldface letters to denote vectors in physical space (i.e.,
$\reals^3$),
and $\| \cdot \|$ to denote the Euclidean norm of vector in the physical space.
Vectors in solution space, such as $q\in \reals^8$, 
where $q$ is the vector of conserved variables for the ideal MHD equations:
$q=\left(\rho,\, \rho \u, \, {\mathcal E}, \B \right)^T$, are not denoted using
boldface letters.}:
\begin{gather}
\label{eqn:MHD}
 \frac{\partial}{\partial t}
  \begin{bmatrix}
    \rho \\ \rho \u \\ \En \\ \B
  \end{bmatrix} +  \nabla \cdot
  \begin{bmatrix} \rho \u \\ \rho \u  \u + \left( {p} + \frac{1}{2}
  \| \B \|^2  \right) {\mathbb I}
    - \B  \B \\
    \u \left(\En + {p} + \frac{1}{2} \| \B \|^2 \right) - 
    	\B \left(\u \cdot \B \right)
    \\ \u  \B - \B  \u
  \end{bmatrix} = 0, \\
  \label{eqn:divfree}
  \nabla \cdot \B = 0,
\end{gather}
where $\rho$, $\rho \u$, and $\En$ are the mass,  momentum,
and energy densities of the plasma system, and $\B$ is the
magnetic field. The thermal pressure, $p$, is related to the
conserved quantities through the ideal gas law:
\begin{equation}
\label{eqn:eos}
      p = \left( \gamma-1 \right) \left( \En - \frac{1}{2} \| \B \|^2 
      	- \frac{1}{2} \rho \| \u \|^2 \right),
\end{equation}
where $\gamma = 5/3$ is the ideal gas constant.

Note that the equation for the magnetic field comes from
Faraday's law:
\begin{equation}
   \B_{,t} + \nabla \times \E = 0,
\end{equation}
where the electric field, $\E$, is approximated
by Ohm's law for a perfect conductor:
\begin{equation}
   \E = \B \times \u.
\end{equation}
Since the electric field is determined entirely from Ohm's law, we 
do not require an evolution equation for it; and thus, the only
other piece that we need from Maxwell's equations is
the divergence-free condition on the magnetic field
\eqref{eqn:divfree}. A complete derivation and discussion
of MHD system \eqref{eqn:MHD}-\eqref{eqn:divfree} can
be found in several standard plasma physics textbooks
(e.g., pages 69--78 of \cite{book:Go98}).

\subsection{Hyperbolicity}
We first note that system \eqref{eqn:MHD}, along with the
equation of state \eqref{eqn:eos}, provides
a full set of equations for the time evolution 
of all eight state variables: $\left( \rho, \, \rho \u, \, \En, \, \B \right)$.
These evolution equations form a hyperbolic system. In particular, the
eigenvalues of the flux Jacobian in some arbitrary
direction ${\bf n}$ ($\| {\bf n} \| = 1$) can be written as follows:
  \begin{gather}
  \lambda^{1,8} = {\bf u} \cdot {\bf n} \mp c_f  \quad \text{(fast magnetosonic),} \qquad
  \lambda^{2,7} = {\bf u} \cdot {\bf n}  \mp c_a  \quad \text{(Alfv$\acute{\text{e}}$n),} \\
  \lambda^{3,6} = {\bf u} \cdot {\bf n}  \mp c_s  \quad \text{(slow magnetosonic),} \qquad
  \lambda^{4} = {\bf u} \cdot {\bf n}  \quad \text{(entropy),} \qquad
  \lambda^{5} = {\bf u} \cdot {\bf n}   \quad \text{(divergence),}
\end{gather}
where 
\begin{gather}
  a \equiv \sqrt{\frac{\gamma p}{\rho}}, \qquad
  c_a \equiv \sqrt{\frac{\left(\B \cdot \n \right)^2}{\rho}}, \qquad
  c_f, c_s \equiv \left\{ \frac{1}{2} \left[ a^2 +
      \frac{\|\B\|^{2}}{\rho} \pm \sqrt{\left(a^2 +
      \frac{\|\B\|^{2}}{\rho} \right)^{2} - 4 a^2
      \frac{\left(\B \cdot \n \right)^2}{\rho}} \right] \right\}^{1/2}.
\end{gather}
The eigenvalues are well-ordered in the sense that
\begin{equation}
  \lambda^{1} \le \lambda^{2} \le \lambda^{3} \le \lambda^{4} \le \lambda^{5} \le \lambda^{6}
  \le \lambda^{7} \le \lambda^{8} \, .
\end{equation}
The fast and slow magnetosonic waves are genuinely nonlinear, while the remaining waves are linearly degenerate.
Note that the so-called {\it divergence-wave} has been
made to travel at the speed $\u \cdot {\bf n}$ via the 
8-wave formulation,
thus restoring Galilean invariance \cite{article:Go72,article:Po94,article:Po99}.
Also note that despite the fact that we use the eigenvalues and eigenvectors
of the 8-wave formulation of the MHD equations, we will still
solve the MHD equations in conservative form (i.e., without the 
``Powell source term''). See \S \ref{subsec:sym_mhd}
and \S \ref{sec:dg_stability} for more discussion on this issue.

\subsection{Symmetric hyperbolic structure of standard hyperbolic conservation laws}
\label{sec:symmhyp}
Consider a conservation
law in $d$ spatial dimensions with $\meq$ number
of conserved variables of the form:
\begin{equation}
\label{eqn:conslaw}
q_{,t} +  \nabla \cdot {\bf F}(q) = 0,
\quad \text{in} \quad \xv \in \Omega \subset \reals^d,
\end{equation}
with appropriate initial and boundary conditions.
In this equation $q: \reals^+ \times \reals^d \mapsto \reals^{\meq}$ 
is the vector of conserved variables and
${\bf F} :\reals^{\meq} \mapsto \reals^{{\meq}\times d}$ is the flux
function. We assume that equation (\ref{eqn:conslaw}) is hyperbolic, 
 meaning that the family of ${\meq} \times {\meq}$ matrices defined by
 \begin{equation}
   \label{eqn:fluxJacobian}
    A(q; \nv) = 
    \frac{\partial \left( \nv \cdot {\bf F} \right)}{\partial q}
\end{equation}
are diagonalizable with real eigenvalues for
all $q$ in the domain of interest and
for all vectors $\nv$ such that $\| \nv \| = 1$.
The matrix $A(q; \nv)$ in \eqref{eqn:fluxJacobian} 
is often referred to as the {\it flux Jacobian}.

Many important hyperbolic conservation laws arising in physical applications
are  equipped with a convex entropy extension that can be used along
with appropriate jump conditions to determine what weak solutions
of \eqref{eqn:conslaw} can be considered admissible. 
Let $\ent: \reals^{\meq} \mapsto \reals$ and $\eflux: \reals^{\meq} \mapsto \reals^d$
denote an entropy function and associated entropy flux that satisfy
the following entropy inequality:
\begin{equation}
\label{eqn:entropy_inequality}
   \ent(q)_{,t} + \nabla \cdot \eflux(q) \le 0,
\end{equation}
where $q$ is a weak solution of \eqref{eqn:conslaw}.
The definitions of $U$ and $\eflux$ are such that inequality \eqref{eqn:entropy_inequality} becomes an exact
equality for classical solutions.

The existence of the above entropy and the associated entropy flux
is closely related to the fact that system \eqref{eqn:conslaw} can
in many cases (e.g., shallow water equations, compressible Euler equations,
and Euler-Maxwell) 
be put in symmetric hyperbolic form via a transformation of the
dependent variables. In particular, we seek a mapping $q(v)$ from $\reals^{\meq}$ to $\reals^{\meq}$ such
that \eqref{eqn:conslaw} becomes
\begin{equation}
\label{eqn:entropy_var_form}
  q_{,v} \, v_{,t} + \sum_{\ell=1}^d F^{\ell}_{,v} \, \, v_{,x^{\ell}} = 0,
\end{equation}
where $q_{,v} \in \reals^{{\meq} \times {\meq}}$ is symmetric positive definite 
and $F^{\ell}_{,v}\in \reals^{{\meq} \times {\meq}}$
  is symmetric for all $\ell = 1, \ldots, d$. 
In several important cases it turns out that the symmetrization
variables, $v$, are directly related to the entropy function
via the relationship:
\begin{equation}
  v = \ent_{,q}.
\end{equation}
For more details on the the symmetrization theory for
nonlinear hyperbolic conservation laws  see for example
\cite{article:God61,article:Moc80,article:Tad87a,article:Ha83b}.

\subsection{Symmetric hyperbolic structure of ideal MHD}
\label{subsec:sym_mhd}
In the case of ideal MHD, a suitable entropy
function is the physical entropy density (divided by $\gamma-1$ for convenience):
\begin{equation}
\label{eqn:entropy_function}
   \ent(q) :=  -\frac{\rho s}{\gamma-1},
\end{equation}
where $s$ is the physical entropy:
\begin{equation}
s:=\log\left(\frac{p}{\rho^\gamma}\right).
\end{equation}
The minus sign in \eqref{eqn:entropy_function} is there to make sure that the entropy function
decreases, which is the usual convention in the theory of hyperbolic conservation laws.
Computing the time derivative of the entropy function and using the
conservation laws from the MHD system \eqref{eqn:MHD} results in the
following equation for the entropy function (for smooth solutions):
\begin{equation}
   \ent(q)_{,t} + \nabla \cdot \eflux(q) + \chi(q)  \left( \nabla \cdot \B \right) = 0,
\end{equation}
where 
\begin{equation}
\label{eqn:chi}
\eflux(q) := \u U \quad \text{and} \quad \chi(q):=\frac{\rho \u \cdot \B}{p}.
\end{equation}
The key point is that the entropy function is only conserved (for smooth solutions)
if the magnetic field is exactly divergence-free: $\nabla \cdot \B = 0$.
Since entropy and symmetrization are connected,
this result suggests that the involution \eqref{eqn:divfree} might play an important
role in the symmetrization of the MHD equations. Indeed it does.

The entropy variables are
\begin{equation}
\label{eqn:entropy_vars}
	v = \ent_{,q} = \left[\frac{\gamma - s}{\gamma -1 } 
	- \frac{\| {\bf u} \|^2}{2T}  , \, \frac{{\bf u}}{T}, \, -\frac{1}{T}, \, \frac{\B}{T} \right]^T,
\end{equation}	
where $T=p/\rho$ is the temperature. A straightforward calculation
shows that using these variables results in an equation of the form
\eqref{eqn:entropy_var_form}, where $q_{,v} \in \reals^{8\times8}$
is indeed symmetric positive definite, but where the matrices
$F^\ell_{,v} \in \reals^{8\times8}$ for $\ell=1,2,3$ are unfortunately
not symmetric. This fact was first observed by Godunov \cite{article:Go72}.

Godunov \cite{article:Go72} provided a remedy for this situation, which
involved adding to \eqref{eqn:entropy_var_form} a term that is proportional
to $\nabla \cdot \B$:
\begin{equation}
\label{eqn:sym_hyp}
\begin{split}
   q_{,v} \, v_{,t} + \sum_{\ell=1}^d F^\ell_{,v} \, v_{,x^\ell} 
   + \chi_{,v} \, \left( \nabla \cdot \B \right) = 0 \quad \Longrightarrow \quad
   q_{,v} \, v_{,t} + \sum_{\ell=1}^d \widetilde{A}^{\ell} \, v_{,x^\ell}  = 0,
   \end{split}
\end{equation}
where $\chi$ is defined by \eqref{eqn:chi} and the matrices $\widetilde{A}^{\ell} \in \reals^{8\times 8}$,
for $\ell=1,2,3$, are symmetric \cite{article:Barth05}.

The additional term in the above expression is sometimes referred to as the 
Powell source term \cite{article:Po94,article:Po99}, since Powell advocated including
this term in numerical simulations of MHD (see
\S\ref{sec:dg_stability} for much more discussion on this point).
The fact that the involution $\nabla \cdot \B = 0$ is needed in order to symmetrize
the MHD equations has direct consequences on the stability of numerical
discretizations of MHD. We will discuss this in more detail in \S \ref{sec:dg_stability}.

We also note the following identities, which will become useful in the DG stability
analysis:
\begin{align}
\label{eqn:entropy_dot1}
   v \cdot q_{,t} &= \ent_{,t}, \\
   \label{eqn:entropy_dot2}
   v \cdot \left( \nabla \cdot {\bf F} \right) &= \nabla \cdot \eflux - \chi \, \nabla \cdot \B, \\
   \label{eqn:entropy_dot3}
   v \cdot \chi_{,v} \, \nabla \cdot \B 
   &=  \chi \,  \nabla \cdot \B.
\end{align}
These identities conspire to produce the following important result in the case of smooth
solutions:
\begin{equation}
 v \cdot \left\{ q_{,t} + \nabla \cdot {\bf F} + \chi_{,v} \, \nabla \cdot {\bf B}  \right\} = \ent_{,t} + \nabla \cdot \eflux = 0.
\end{equation}
In words, this result says that multiplying the augmented MHD equations (i.e., MHD with the 
``Powell source term'')
 by the symmetrization variables yields the conservation law for the
 entropy function.
 Finally, we introduce the dual entropy and dual entropy flux:
 \begin{alignat}{2}
  \label{eqn:dual_entropy}
  \ent + \dualent &= v \cdot q &\qquad \Longrightarrow \qquad 
  \dualent &:= \rho + \frac{\| \B \|^2}{2T},  \\
  \label{eqn:dual_entropy_flux}
  \eflux + \dualeflux &= v \cdot {\bf F} + \chi \, {\bf B} 
  &\qquad \Longrightarrow \qquad \dualeflux &:=  \rho \u + \frac{\u \| \B \|^2}{2T},
 \end{alignat}
 which satisfy the following relationships:
 \begin{equation}
 \label{eqn:dual_vars}
  q = \dualent_{,v} \qquad \text{and} \qquad
  {\bf F} = \dualeflux_{,v} - \chi_{,v} \, {\bf B}.
 \end{equation}

\section{Discontinuous Galerkin finite element methods}
\label{sec:dg_method}
The modern form of the discontinuous Galerkin (DG) finite element method
for solving hyperbolic PDEs of the form \eqref{eqn:conslaw}
was developed in a series of papers by Bernardo Cockburn, 
Chi-Wang Shu, and their collaborators (see e.g., \cite{article:CoShu5}). 
In this section we briefly
review the general framework of the DG method in \S\ref{sec:dg_general}, and give
further details for 2D Cartesian meshes in \S\ref{sec:dg-cart-2d} and 3D Cartesian meshes
in  \S\ref{sec:dg-cart-3d}. This section also serves to explain the notation
used throughout this paper.

\subsection{General framework}
\label{sec:dg_general}
Let $\Omega \subset \reals^d$ be a polygonal domain with boundary $\partial \Omega$. 
We discretize $\Omega$ using a finite set of non-overlapping elements, $\Tm_i$, such
that  $\Omega = \cup_{i=1}^N \Tm_i$. 
Let ${P}^{\, \deg}\left(\reals^d\right)$ denote the set of polynomials from $\reals^d$ to $\reals$ 
with maximal polynomial degree $\deg$; for example, the polynomial
$x^a y^b z^c$, where $a,b,c$ are all 
non-negative integers, is in ${P}^{\, \deg}\left(\reals^3\right)$ if and only if $a+b+c \le \deg$.
On the mesh of $N$ elements we define the {\it broken} finite element space
\begin{equation}
\label{eqn:broken_space}
    \WS^{h,M_{\text{eq}}}_{\, \deg} = \left\{ w^h \in \left[ L^{\infty}(\Omega) \right]^{M_{\text{eq}}}: \,
    w^h |_{\Tm_i} \in \left[ {P}^{\, \deg} \right]^{M_{\text{eq}}}, \, \forall \Tm_i \in \Tm^h \right\},
\end{equation}
where $h$ is the grid spacing, $M_{\text{eq}}$ is the number of solution variables, and $\deg$ is the 
maximal polynomial degree in the finite element representation.
The above expression means that $w^h \in \WS^{h,M_{\text{eq}}}_{\deg}$ has
$M_{\text{eq}}$ components, each of which when restricted to some element $\Tm_i$
is a polynomial of degree at most $\deg$ and no continuity is assumed
across element faces (or edges in 2D). 

The semi-discrete DG finite element method
is obtained by multiplying \eqref{eqn:conslaw} by a test function
$w^h({\bf x}) \in \WS^{h,\meq}_{\deg}$, replacing the exact solution $q(t,{\bf x})$ by
the trial function $q^h(t,{\bf x}) \in  \WS^{h,\meq}_{\deg}$, and carrying out the appropriate integrations-by-part. 
The result of this is the following semi-discrete variational problem
for all $t\in[0,T]$:
\begin{equation}
\label{eqn:semidiscrete_fem1}
\begin{split}
     \text{Find} \quad &q^h \in \WS^{h,\meq}_{\deg} \quad \text{such that} \\
     &B_{\text{DG}}\left(q^h,w^h \right) = 0 \qquad
       \text{for all} \quad w^h \in \WS^{h,\meq}_{\deg} \quad \text{and} \quad t\in[0,T],
       \end{split}
\end{equation}
where
\begin{equation}
\begin{split}
   B_{\text{DG}}\left(q^h,w^h \right) :=  &\sum_{i=1}^N \ \iiint_{\Tm_i}  \left\{ w^h \cdot
q^h_{,t} - \nabla w^h : {\bf F}\left(q^h \right) \right\} \, d{\bf x} 
 +   \sum_{i=1}^N  \varoiint_{\partial \Tm_i}
       w^h_{-}  \cdot \NF\left(q^h_{-}({\bf s}), \, q^h_{+}({\bf s}); \, {\bf n}\right) \, 
       d{\bf s}.
\label{eqn:semidiscrete_fem2}
\end{split}
\end{equation}
In the above expression, ${\bf n}$ is a unit normal that is outward pointing
relative to element $\Tm_i$, and $q^h_{-}({\bf s})$ and $q^h_{+}({\bf s})$ denote the approximate
solutions on either side of the boundary $\partial \Tm_i$, where ${\bf s}$ are
are coordinates on $\partial \Tm_i$ and
the subscript minus (plus) sign means {\it interior to} ({\it exterior to}) element $\Tm_i$. 
Similarly, $w^h_{-}$ represents
the test function evaluated on  the boundary $\partial \Tm_i$ on the
side that is interior to $\Tm_i$.
In this expression $\NF\left(q^h_-, \, q^h_+; \, {\bf n}\right) : 
\reals^m \times \reals^m \times \reals^d \mapsto \reals^m$ 
denotes a numerical flux function with the following two properties:
\begin{itemize}
  \item Consistency:  
  \[ 
  \NF\left(q, \, q; \, {\bf n} \right) = {\bf F}(q) \cdot {\bf n},
  \]
  \item Conservation: 
  \[
  \NF\left(q^h_-, \, q^h_+; \, {\bf n} \right) = -\NF\left(q^h_+, \, q^h_-; \, -{\bf n} \right).
  \]
\end{itemize}
As a matter of practice we use in this work the local Lax-Friedrichs (LLF) flux \cite{article:Ru61}:
\begin{equation}
\label{eqn:LLF}
\NF\left(q^h_-, \, q^h_+; \, {\bf n} \right) = \frac{1}{2} \left[ {\bf n} \cdot \left( {\bf F}(q^h_+) + {\bf F}(q^h_-) \right)
- \alpha \left( q^h_+ - q^h_- \right) \right],
\end{equation}
where $\alpha$ is an estimate of the maximum local wave speed.
In the next two subsections we give details for the above described numerical
scheme on 2D (\S \ref{sec:dg-cart-2d}) and
3D (\S \ref{sec:dg-cart-3d}) Cartesian meshes.


\subsection{2D Cartesian meshes}
\label{sec:dg-cart-2d}
Let $\Tm^h$ be a Cartesian grid over $\Omega = 
[a_x, \, b_x] \times [a_y, \, b_y]$,
with uniform grid spacing $\Delta x$ and $\Delta y$ in each
coordinate direction, where $h=\max\left(\Delta x, \Delta y\right)$.
The mesh elements are centered 
at the coordinates
\begin{equation}
          x_i =  a_x  +  \left( i - \frac{1}{2} \right)  \Delta x
          \quad \text{and} \quad
          y_j =  a_y  +  \left( j - \frac{1}{2} \right) \Delta y.
\end{equation}
Each element can be mapped to the canonical element $(\xi, \eta) \in [-1,1] \times [-1,1]$ via the linear transformation:
\begin{equation}
	x = x_i + \xi \, \frac{\Delta x}{2}, \quad 
	y = y_j + \eta \, \frac{\Delta y}{2}.
\end{equation}
We define an orthonormal set of polynomial basis functions
that span the broken finite element space
\eqref{eqn:broken_space}. In particular, up to degree three
these basis functions can be written in canonical coordinates as
\begin{equation}
\label{eqn:test-fun-cart}
\begin{split}
\varphi^{(\ell)} \in \Biggl\{ &1, \, \, \,  \sqrt{3} \, \xi, \, \, \,  \sqrt{3} \, \eta, 
  \, \, \,   3 \, \xi \eta, \, \, \,  \frac{\sqrt{5}}{2} \left( 3 \xi^2 - 1 \right),
\, \, \,  \frac{\sqrt{5}}{2} \left( 3 \eta^2 - 1 \right), \\ 
&
\frac{\sqrt{7}}{2} (5 \xi^3 - 3 \xi), \, \, \,
\frac{\sqrt{15}}{2} \eta \, (3 \xi^2 - 1), \, \, \,
\frac{\sqrt{15}}{2} \xi \, (3 \eta^2 - 1), \, \, \,
\frac{\sqrt{7}}{2} (5 \eta^3 - 3 \eta) \Biggr\}.
\end{split}
\end{equation}
These basis functions have been orthonormalized with respect
to the following inner product:
\begin{equation}
   \Bigl\langle \varphi^{(m)}, \, \varphi^{(n)} \Bigr\rangle :=  
   \frac{1}{4}  \int_{-1}^{1} \int_{-1}^{1}
    \varphi^{(m)}(\xi,\eta) \, \varphi^{(n)}(\xi,\eta) \, d\xi \, d\eta = 
     \delta_{mn}.
\end{equation}

We  look for approximate solutions of (\ref{eqn:conslaw})
that have the following form:
\begin{equation}
\label{eqn:q_ansatz}
q^h(t, x(\xi), y(\eta)) \Bigl|_{\Tm_{ij}}  :=
  \sum_{m=1}^{M(M+1)/2}  Q^{(m)}_{ij}(t) \, \varphi^{(m)}(\xi, \eta),
\end{equation}
where $M$ is the desired order of accuracy in space and $Q^{(\ell)}_{ij}(t)$
represents the $M^{\text{th}}$ Legendre coefficient.
The Legendre coefficients of the initial conditions at $t=0$ are 
determined from the $L^2$-projection of $q_0(x,y)$ onto the Legendre basis functions:
\begin{align}
\label{eqn:l2project}
Q^{(\ell)}_{ij}(0)  &:= \Bigl\langle q_0\left(x_i + 0.5 \, {\Delta x} \, \xi, \, y_j + 0.5 \, {\Delta y} \, \eta \right), \, \varphi^{(\ell)}(\xi,\eta) \Bigr\rangle.
\end{align}
In practice, these double integrals are evaluated using
standard 2D Gaussian quadrature rules involving $M^2$ points. 

In order to determine the Legendre coeficients for $t>0$,
we take the semi-discrete equation
\eqref{eqn:semidiscrete_fem1}--\eqref{eqn:semidiscrete_fem2}
and replace the trial function, $q^h$, by \eqref{eqn:q_ansatz} and the test function,
$w^h$, by Legendre basis function $ \varphi^{(\ell)}$ that are supported on a single element.
This results in the following semi-discrete DG method:
\begin{equation}
\label{eqn:semidiscrete}
\frac{d}{dt} \,
   Q^{(\ell)}_{ij} = {\mathcal L}^{(\ell)}_{ij} := N^{(\ell)}_{ij}
- \frac{\Delta {F}_{ij}^{(\ell)}}{\Delta x}
- \frac{\Delta {G}_{ij}^{(\ell)}}{\Delta y} ,
\end{equation}
where 
\begin{gather}
	\label{eqn:Nvals}
	N^{(\ell)}_{ij} := \frac{1}{2}
   \int_{-1}^{1} \int_{-1}^{1}
   \left[ \, \frac{1}{\Delta x} \, \varphi^{(\ell)}_{, \xi} \, F^1(q^h) 
   	+ \frac{1}{\Delta y} \, \varphi^{(\ell)}_{, \eta} \,  F^2(q^h) \, \right] \, 
   d\xi  \, 
   d\eta, \\
	\label{eqn:Fl1}
	\Delta {F}^{(\ell)}_{ij} := 
	\left[ \frac{1}{2} \int_{-1}^{1}
	\varphi^{(\ell)} \,
		 \NF\left(q^h_{-}, \, q^h_{+}; \, {\bf e}^1 \right)   \, 
		 d\eta \right]_{\xi=-1}^{\xi=1},  \quad
		\Delta {G}^{(\ell)}_{ij} := \left[
		 \frac{1}{2} \int_{-1}^{1}
	\varphi^{(\ell)} \,
		 \NF\left(q^h_{-}, \, q^h_{+}; \, {\bf e}^2 \right)   \, d\xi
		 \right]_{\eta=-1}^{\eta=1},
\end{gather}
and
\begin{equation}
   {\bf e}^1 := \bigl( 1, \, 0 \bigr)^T \qquad \text{and} \qquad
   {\bf e}^2 := \bigl( 0, \, 1 \bigr)^T.
\end{equation}
The integrals in (\ref{eqn:Nvals}) can
be numerically approximated via standard 
2D Gaussian quadrature rules involving $(M-1)^2$ points. 
The integrals in \eqref{eqn:Fl1} can
be approximated with the standard 1D Gauss quadrature
rule  involving $M$ points, and the numerical fluxes
in these expressions are evaluated at each Gaussian
quadrature point using the local Lax-Friedrichs (LLF) flux given by
\eqref{eqn:LLF}.

\subsection{3D Cartesian meshes}
\label{sec:dg-cart-3d}
Let $\Tm^h$ be a Cartesian grid over $\Omega = 
[a_x, \, b_x] \times [a_y, \, b_y] \times [a_z, \, b_z]$,
with uniform grid spacing $\Delta x$, $\Delta y$, and $\Delta z$ in each
coordinate direction, where $h=\max\left(\Delta x, \Delta y, \Delta z\right)$.
The mesh elements are centered 
at the coordinates
\begin{equation}
          x_i =  a_x  +  \left( i - \frac{1}{2} \right)  \Delta x, \quad
          y_j =  a_y  +  \left( j - \frac{1}{2} \right)  \Delta y,
          \quad \text{and} \quad
          z_k =  a_z  +  \left( k - \frac{1}{2} \right) \Delta z.
\end{equation}
Each element can be mapped to the canonical element $(\xi, \eta, \zeta) \in [-1,1] \times [-1,1] \times [-1,1]$ 
via the linear transformation:
\begin{equation}
	x = x_i + \xi \, \frac{\Delta x}{2}, \quad 
	y = y_j + \eta \, \frac{\Delta y}{2}, \quad
	z = z_k + \zeta \, \frac{\Delta y}{2}.
\end{equation}
We define an orthonormal set of polynomial basis functions
that span the broken finite element space
\eqref{eqn:broken_space}. In particular, up to degree three
these basis functions can be written in canonical coordinates as
\begin{equation}
\label{eqn:test-fun-cart-3d}
\begin{split}
\varphi^{(\ell)} \in \Biggl\{ 
&1, 
\, \, \,  \sqrt{3} \, \xi, 
\, \, \,  \sqrt{3} \, \eta, 
\, \, \,  \sqrt{3} \, \zeta, 
\, \, \,   3 \, \xi \eta,   
\, \, \,   3 \, \xi \zeta,   
\, \, \,   3 \, \eta \zeta, 
\, \, \,  \frac{\sqrt{5}}{2} \left( 3 \xi^2 - 1 \right),
\, \, \,  \frac{\sqrt{5}}{2} \left( 3 \eta^2 - 1 \right),
\, \, \,  \frac{\sqrt{5}}{2} \left( 3 \zeta^2 - 1 \right), \\ 
&
\frac{\sqrt{15}}{2} \eta \left( 3 \xi^2 - 1 \right), \, \, \, 
\frac{\sqrt{15}}{2} \zeta \left( 3 \xi^2 - 1 \right), \, \, \, 
\frac{\sqrt{15}}{2} \xi \left( 3 \eta^2 - 1 \right), \, \, \, 
\frac{\sqrt{15}}{2} \zeta \left( 3 \eta^2 - 1 \right), \, \, \, 
\frac{\sqrt{15}}{2} \xi \left( 3 \zeta^2 - 1 \right), \\
&
\frac{\sqrt{15}}{2} \eta \left( 3 \zeta^2 - 1 \right), \, \, \, 
3 \xi  \eta \zeta, \, \, \, 
\frac{\sqrt{7}}{2} (5 \xi^3 - 3 \xi), \, \, \,
\frac{\sqrt{7}}{2} (5 \eta^3 - 3 \eta), \, \, \, 
\frac{\sqrt{7}}{2} (5 \zeta^3 - 3 \zeta) \Biggr\}.
\end{split}
\end{equation}
These basis functions have been orthonormalized with respect
to the following inner product:
\begin{equation}
   \Bigl\langle \varphi^{(m)}, \, \varphi^{(n)} \Bigr\rangle :=  
   \frac{1}{8}  \int_{-1}^{1} \int_{-1}^{1}  \int_{-1}^{1}
    \varphi^{(m)}(\xi,\eta,\zeta) \, \varphi^{(n)}(\xi,\eta,\zeta) \, d\xi \, d\eta \, d\zeta = 
     \delta_{mn}.
\end{equation}

We  look for approximate solutions of (\ref{eqn:conslaw})
that have the following form:
\begin{equation}
\label{eqn:q_ansatz-3d}
q^h(t, x(\xi), y(\eta), z(\zeta)) \Bigl|_{\Tm_{ijk}}  :=
  \sum_{m=1}^{M(M+1)(M+2)/6}  Q^{(m)}_{ijk}(t) \, \varphi^{(m)}(\xi, \eta, \zeta),
\end{equation}
where $M$ is the desired order of accuracy in space and $Q^{(\ell)}_{ijk}(t)$
represents the $m^{\text{th}}$ Legendre coefficient.
The Legendre coefficients of the initial conditions at $t=0$ are 
determined from the $L^2$-projection of $q_0(x,y,z)$ onto the Legendre basis functions:
\begin{align}
\label{eqn:l2project-3d}
Q^{(\ell)}_{ijk}(0)  &:= \Bigl\langle q_0\left(x_i + 0.5 \, {\Delta x} \, \xi, \, y_j + 0.5 \, {\Delta y} \, \eta,
\, z_k + 0.5 \, {\Delta z} \, \zeta \right), \, \varphi^{(\ell)}(\xi,\eta,\zeta) \Bigr\rangle.
\end{align}
In practice, these triple integrals are evaluated using
standard 3D Gaussian quadrature rules involving $M^3$ points. 

we take the semi-discrete equation
\eqref{eqn:semidiscrete_fem1}--\eqref{eqn:semidiscrete_fem2}
and replace the trial function, $q^h$, by \eqref{eqn:q_ansatz-3d} and the test function,
$w^h$, by Legendre basis function $ \varphi^{(\ell)}$ that are supported on a single element.
This results in the following semi-discrete DG method:
\begin{equation}
\label{eqn:semidiscrete-3d}
\frac{d}{dt} \,
   Q^{(\ell)}_{ijk} = {\mathcal L}^{(\ell)}_{ijk} := N^{(\ell)}_{ijk}
- \frac{\Delta {F}_{ijk}^{(\ell)}}{\Delta x}
- \frac{\Delta {G}_{ijk}^{(\ell)}}{\Delta y} 
- \frac{\Delta {H}_{ijk}^{(\ell)}}{\Delta z},
\end{equation}
where 
\begin{align}
	\label{eqn:Nvals-3d}
	N^{(\ell)}_{ijk} &:= \frac{1}{4}
   \int_{-1}^{1} \int_{-1}^{1} \int_{-1}^{1}
   \left[ \, \frac{1}{\Delta x} \, \varphi^{(\ell)}_{, \xi} \, F^1(q^h) 
   	+ \frac{1}{\Delta y} \, \varphi^{(\ell)}_{, \eta} \,  F^2(q^h)
	+ \frac{1}{\Delta z} \, \varphi^{(\ell)}_{, \zeta} \,  F^3(q^h) \, \right] \, 
   d\xi  \, 
   d\eta \, d\zeta, \\
	\label{eqn:Fl1-3d}
	\Delta {F}^{(\ell)}_{ijk} &:= 
	\left[ \frac{1}{4} \int_{-1}^{1} \int_{-1}^{1}
	\varphi^{(\ell)} \,
		 \NF\left(q^h_{-}, \, q^h_{+}; \, {\bf e}^1 \right)   \, 
		 d\eta \, d\zeta \right]_{\xi=-1}^{\xi=1},  \\ 
	\label{eqn:Gl1-3d}
		\Delta {G}^{(\ell)}_{ijk} &:= \left[
		 \frac{1}{4} \int_{-1}^{1} \int_{-1}^{1}
	\varphi^{(\ell)} \,
		 \NF\left(q^h_{-}, \, q^h_{+}; \, {\bf e}^2 \right)   \, d\xi \, d\zeta
		 \right]_{\eta=-1}^{\eta=1}, \\
	\label{eqn:Hl1-3d}
		\Delta {H}^{(\ell)}_{ijk} &:= \left[
		 \frac{1}{4} \int_{-1}^{1} \int_{-1}^{1}
	\varphi^{(\ell)} \,
		 \NF\left(q^h_{-}, \, q^h_{+}; \, {\bf e}^3 \right)   \, d\xi \, d\eta
		 \right]_{\zeta=-1}^{\zeta=1},
\end{align}
and
\begin{equation}
   {\bf e}^1 := \bigl( 1, \, 0, \, 0 \bigr)^T, \qquad
   {\bf e}^2 := \bigl( 0, \, 1, \, 0 \bigr)^T, \qquad \text{and} \qquad
   {\bf e}^3 := \bigl( 0, \, 0, \, 1 \bigr)^T.
\end{equation}
The integrals in (\ref{eqn:Nvals-3d}) can
be numerically approximated via standard 
3D Gaussian quadrature rules involving $(M-1)^3$ points. 
The integrals in \eqref{eqn:Fl1-3d}, \eqref{eqn:Gl1-3d}, and
\eqref{eqn:Hl1-3d} can be approximated with the standard 2D Gauss quadrature
rule  involving $M^2$ points, and the numerical fluxes
in these expressions are evaluated at each Gaussian
quadrature point using the local Lax-Friedrichs (LLF) flux given by
\eqref{eqn:LLF}.

\subsection{Time-stepping}
\label{sec:timestepping}
The time-stepping is handled via standard total-variation diminishing
Runge-Kutta  (TVD-RK) methods \cite{article:GoShu98,gottliebShuTadmor01}. 
In particular, in this work we make use of the third order accurate version:
\begin{align}
Q^{\star} &= Q^n + \Delta t \, {\mathcal L}({Q}^n), \quad
   Q^{\star \star} = \frac{3}{4} Q^n + \frac{1}{4} Q^{\star}
          + \frac{1}{4} \Delta t \, {\mathcal L}({Q}^{\star}), \quad
   Q^{n+1} = \frac{1}{3} Q^n + \frac{2}{3} Q^{\star \star}
          + \frac{2}{3} \Delta t \, {\mathcal L}({Q}^{\star \star}).
\end{align}

\subsection{Limiters}
\label{sec:limiters}
High-order methods on nonlinear hyperbolic conservation
laws generally require the use of additional limiters
in order to stabilize them on problems with shocks.
In this section we briefly describe the limiting strategies used in this
work.  In the formulas described below we will
make use of the minmod function:
\begin{equation}
  \text{mm}(a, \, b, \, c) = 
  \begin{cases}
  	\text{sgn}(a) \, \text{min}\left( |a|, \, |b|, \, |c| \right) &
	\text{if} \quad \text{sgn}(a) = \text{sgn}(b) = \text{sgn}(c), \\
	0 & \text{otherwise},
  \end{cases}
\end{equation}
where the $\text{sgn}(a)=\pm 1$ is the sign of the real number $a$.

\subsubsection{2D Cartesian meshes}
On 2D Cartesian meshes we make use of moment limiters similar to those
introduced by Krivodonova \cite{article:Kriv07} in order to stabilize the DG 
method in the vicinity of shocks. In the formulas below we will
make use of the following matrices:
\begin{align}
R^{x}_{ij} := R^{x}\left( Q^{(1)}_{ij} \right) \quad \text{and} \quad
R^{y}_{ij} := R^{y}\left( Q^{(1)}_{ij} \right), 
\end{align}
where $R^{x}$ and $R^{y}$ are the matrices of right eigenvectors of $F_{,q} \cdot {\bf e}^1$ and
$F_{,q} \cdot {\bf e}^2$, respectively. The corresponding matrices of left eigenvectors are denoted as follows:
\begin{equation}
 L^{x}_{ij} := \bigl( R^{x}_{ij} \bigr)^{-1} \quad \text{and} \quad
 L^{y}_{ij} := \bigl( R^{y}_{ij} \bigr)^{-1}.
\end{equation}
The specific matrices used in this work 
are those introduced by Barth \cite{article:barth99}, who
developed right-eigenvector scalings for ideal MHD that have optimal direction-independent matrix norms.

The procedure we advocate
can be summarized as follows:
\begin{enumerate}
\item Compute the following conversions from conservative to characteristic variables:
\begin{alignat}{4}
	U^{(1)}_{ij} &:=& \, L^{x}_{ij} \, Q^{(2)}_{ij},  
        \quad \Delta_{+} U^{(1)}_{ij} &:=&  \, L^{x}_{ij} \, \left( {Q^{(1)}_{i+1 \, j}-Q^{(1)}_{ij}} \right), 
        \quad  \Delta_{-} U^{(1)}_{ij} &:=&  \, L^{x}_{ij} \, \left( {Q^{(1)}_{ij}-Q^{(1)}_{i-1 \, j}} \right), \\
       U^{(2)}_{ij} &:=& L^{x}_{ij} \, Q^{(4)}_{ij},  
        \quad \Delta_{+} U^{(2)}_{ij} &:=&  \, L^{x}_{ij} \, \left( {Q^{(2)}_{i \, j+1}-Q^{(2)}_{ij}} \right), 
        \quad  \Delta_{-} U^{(2)}_{ij} &:=&  \, L^{x}_{ij} \, \left( {Q^{(2)}_{ij}-Q^{(2)}_{i \, j-1}} \right),  \\
        U^{(3)}_{ij} &:=& L^{x}_{ij} \, Q^{(5)}_{ij},  
        \quad \Delta_{+} U^{(3)}_{ij} &:=&  \, L^{x}_{ij} \, \left( {Q^{(2)}_{i+1 \, j}-Q^{(2)}_{ij}} \right), 
        \quad  \Delta_{-} U^{(3)}_{ij} &:=&  \, L^{x}_{ij} \, \left( {Q^{(2)}_{ij}-Q^{(2)}_{i-1 \, j}} \right),  \\  
        V^{(1)}_{ij} &:=& L^{y}_{ij} \, Q^{(3)}_{ij},  
        \quad \Delta_{+} V^{(1)}_{ij} &:=&  \, L^{y}_{ij} \, \left( {Q^{(1)}_{i \, j+1}-Q^{(1)}_{ij}} \right), 
        \quad  \Delta_{-} V^{(1)}_{ij} &:=&  \, L^{y}_{ij} \, \left( {Q^{(1)}_{ij}-Q^{(1)}_{i \, j-1}} \right), \\
        V^{(2)}_{ij} &:=& L^{y}_{ij} \, Q^{(4)}_{ij},  
        \quad \Delta_{+} V^{(2)}_{ij} &:=&  \, L^{y}_{ij} \, \left( {Q^{(3)}_{i+1 \, j}-Q^{(3)}_{ij}} \right), 
        \quad  \Delta_{-} V^{(2)}_{ij} &:=&  \, L^{y}_{ij} \, \left( {Q^{(3)}_{ij}-Q^{(3)}_{i-1 \, j}} \right),  \\
        V^{(3)}_{ij} &:=& L^{y}_{ij} \, Q^{(6)}_{ij},  
        \quad \Delta_{+} V^{(3)}_{ij} &:=&  \, L^{y}_{ij} \, \left( {Q^{(3)}_{i \, j+1}-Q^{(3)}_{ij}} \right), 
        \quad  \Delta_{-} V^{(3)}_{ij} &:=&  \, L^{y}_{ij} \, \left( {Q^{(3)}_{ij}-Q^{(3)}_{i \, j-1}} \right).      
\end{alignat}
\item Let $u^{(\ell)}_{ij}$ be some component of $U^{(\ell)}_{ij}$. Similarly,
	let $\Delta_{-} u^{(\ell)}_{ij}$, $\Delta_{+} u^{(\ell)}_{ij}$, $v^{(\ell)}_{ij}$,
	$\Delta_{-} v^{(\ell)}_{ij}$, and $\Delta_{+} v^{(\ell)}_{ij}$ be 
	the corresponding components of $\Delta_{-} U^{(\ell)}_{ij}$, $\Delta_{+} U^{(\ell)}_{ij}$,
	$V^{(\ell)}_{ij}$, $\Delta_{-} V^{(\ell)}_{ij}$, and $\Delta_{+} V^{(\ell)}_{ij}$.
\item Loop over all components and limit as follows:
\begin{algorithmic}
\STATE (a) Limit quadratic terms:
\begin{align}
     \tilde{u}^{(3)}_{ij} &:= \text{mm} \left\{ u^{(3)}_{ij}, \, \alpha_3 \,  \Delta_{-} u^{(3)}_{ij}, \,
     \alpha_3 \,  \Delta_{+} u^{(3)}_{ij}  \right\} \quad \text{and} \quad
     \tilde{v}^{(3)}_{ij} := \text{mm} \left\{ v^{(3)}_{ij}, \, \alpha_3 \,  \Delta_{-} v^{(3)}_{ij}, \,
     \alpha_3 \,  \Delta_{+} v^{(3)}_{ij}  \right\}.
\end{align}
\STATE (b) If \, $\left|  \tilde{u}^{(3)}_{ij} - u^{(3)}_{ij} \right| + \left| \tilde{v}^{(3)}_{ij} - v^{(3)}_{ij} \right| > 0$  \,
		or \, $\left|  {u}^{(3)}_{ij} \right| = 0$ \, or \, $\left|  {v}^{(3)}_{ij} \right| = 0$ \,
	 then limit all the other terms ($\ell=1,2$):
\begin{align}
     \tilde{u}^{(\ell)}_{ij} &:= \text{mm} \left\{ u^{(\ell)}_{ij}, \, \alpha_{\ell} \,  \Delta_{-} u^{(\ell)}_{ij}, \,
     \alpha_{\ell} \,  \Delta_{+} u^{(\ell)}_{ij}  \right\} \quad \text{and} \quad
     \tilde{v}^{(\ell)}_{ij} := \text{mm} \left\{ v^{(\ell)}_{ij}, \, \alpha_{\ell} \,  \Delta_{-} v^{(\ell)}_{ij}, \,
     \alpha_{\ell} \,  \Delta_{+} v^{(\ell)}_{ij}  \right\},
\end{align}
\STATE where in this work we take for $\alpha_{\ell}$ numbers in the range
	suggested by Krivodonova \cite{article:Kriv07}:  $\alpha_3 = \alpha_2 = \sqrt{\frac{3}{20}}$ \, and \, $\alpha_1 = \sqrt{\frac{1}{12}}$.
\end{algorithmic}
\item Convert characteristic values back to conservative values:
\begin{equation}
\begin{split}
     {Q}^{(2)}_{ij} \leftarrow R^{x}_{ij} \, \tilde{U}^{(1)}_{ij}, \quad & {Q}^{(3)}_{ij} \leftarrow R^{y}_{ij} \, \tilde{V}^{(1)}_{ij}, \quad
     {Q}^{(4)}_{ij} \leftarrow \text{mm} \left\{ R^{x}_{ij} \, \tilde{U}^{(2)}_{ij}, \, R^{y}_{ij} \, \tilde{V}^{(2)}_{ij} \right\}, \\
     & {Q}^{(5)}_{ij} \leftarrow R^{x}_{ij} \, \tilde{U}^{(3)}_{ij}, \quad
     {Q}^{(6)}_{ij} \leftarrow R^{y}_{ij} \, \tilde{V}^{(3)}_{ij}.
\end{split}
\end{equation}
\end{enumerate}

\subsubsection{3D Cartesian meshes}
In the formulas below we will
make use of the following matrices:
\begin{align}
R^{x}_{ijk} := R^{x}_{ijk}\left( Q^{(1)}_{ijk} \right), \quad
R^{y}_{ijk} := R^{y}_{ijk}\left( Q^{(1)}_{ijk} \right), \quad \text{and} \quad
R^{z}_{ijk} := R^{z}_{ijk}\left( Q^{(1)}_{ijk} \right),
\end{align}
as well as the corresponding matrices of left eigenvectors:
\begin{equation}
 L^{x}_{ijk} := \bigl( R^{x}_{ijk} \bigr)^{-1}, \quad
 L^{y}_{ijk} := \bigl( R^{y}_{ijk} \bigr)^{-1}, \quad \text{and} \quad
 L^{z}_{ijk} := \bigl( R^{z}_{ijk} \bigr)^{-1}.
\end{equation}
The specific left and right-eigenvector scalings used in this work are those
developed by Barth \cite{article:barth99}.
The procedure we advocate
can be summarized as follows:
\begin{enumerate}
\item Compute the following conversions from conservative to characteristic variables:
\begin{alignat}{4}
        U^{(1)}_{ijk} &:=& \, L^{x}_{ijk} \, Q^{(2)}_{ijk},  
        \quad \Delta_{+} U^{(1)}_{ijk} &:=& \, L^{x}_{ijk} \, \left( {Q^{(1)}_{i+1 \, j k}-Q^{(1)}_{ijk}} \right), 
        \quad  \Delta_{-} U^{(1)}_{ijk} &:=& \, L^{x}_{ijk} \, \left( {Q^{(1)}_{ijk}-Q^{(1)}_{i-1 \, j k}} \right), \\
       U^{(2)}_{ijk} &:=& L^{x}_{ijk} \, Q^{(5)}_{ijk},  
        \quad \Delta_{+} U^{(2)}_{ijk} &:=& \, L^{x}_{ijk} \, \left( {Q^{(2)}_{i \, j+1 \, k}-Q^{(2)}_{ijk}} \right), 
        \quad  \Delta_{-} U^{(2)}_{ijk} &:=& \, L^{x}_{ijk} \, \left( {Q^{(2)}_{ijk}-Q^{(2)}_{i \, j-1 \, k}} \right),  \\
        U^{(3)}_{ijk} &:=& \, L^{x}_{ijk} \, Q^{(6)}_{ijk},  
        \quad \Delta_{+} U^{(3)}_{ijk} &:=& \, L^{x}_{ijk} \, \left( {Q^{(2)}_{i \, j k+1}-Q^{(2)}_{ijk}} \right), 
        \quad  \Delta_{-} U^{(3)}_{ijk} &:=& \, L^{x}_{ijk} \, \left( {Q^{(2)}_{ijk}-Q^{(2)}_{i j k-1}} \right),  \\  
        U^{(4)}_{ijk} &:=& \, L^{x}_{ijk} \, Q^{(8)}_{ijk},  
        \quad \Delta_{+} U^{(4)}_{ijk} &:=& \, L^{x}_{ijk} \, \left( {Q^{(2)}_{i+1 \, j k}-Q^{(2)}_{ijk}} \right), 
        \quad  \Delta_{-} U^{(4)}_{ijk} &:=& \, L^{x}_{ijk} \, \left( {Q^{(2)}_{ijk}-Q^{(2)}_{i-1 \, j k}} \right),  \\  
	V^{(1)}_{ijk} &:=& \, L^{y}_{ijk} \, Q^{(3)}_{ijk},  
        \quad \Delta_{+} V^{(1)}_{ijk} &:=& \, L^{y}_{ijk} \, \left( {Q^{(1)}_{i \, j+1 \, k}-Q^{(1)}_{ijk}} \right), 
        \quad  \Delta_{-} V^{(1)}_{ijk} &:=& \, L^{y}_{ijk} \, \left( {Q^{(1)}_{ijk}-Q^{(1)}_{i \, j-1 \, k}} \right), \\
       V^{(2)}_{ijk} &:=& \, L^{y}_{ijk} \, Q^{(5)}_{ijk},  
        \quad \Delta_{+} V^{(2)}_{ijk} &:=& \, L^{y}_{ijk} \, \left( {Q^{(3)}_{i+1 \, j \, k}-Q^{(3)}_{ijk}} \right), 
        \quad  \Delta_{-} V^{(2)}_{ijk} &:=& \, L^{y}_{ijk} \, \left( {Q^{(3)}_{ijk}-Q^{(3)}_{i-1 \, j k}} \right),  \\
        V^{(3)}_{ijk} &:=& \, L^{y}_{ijk} \, Q^{(7)}_{ijk},  
        \quad \Delta_{+} V^{(3)}_{ijk} &:=& \, L^{y}_{ijk} \, \left( {Q^{(3)}_{i j \, k+1}-Q^{(3)}_{ijk}} \right), 
        \quad  \Delta_{-} V^{(3)}_{ijk} &:=& \, L^{y}_{ijk} \, \left( {Q^{(3)}_{ijk}-Q^{(3)}_{i j \, k-1}} \right),  \\  
        V^{(4)}_{ijk} &:=& \, L^{y}_{ijk} \, Q^{(9)}_{ijk},  
        \quad \Delta_{+} V^{(4)}_{ijk} &:=& \, L^{y}_{ijk} \, \left( {Q^{(3)}_{i \, j+1 \, k}-Q^{(3)}_{ijk}} \right), 
        \quad  \Delta_{-} V^{(4)}_{ijk} &:=& \, L^{y}_{ijk} \, \left( {Q^{(3)}_{ijk}-Q^{(3)}_{i \, j-1 \, k}} \right),  \\  
        W^{(1)}_{ijk} &:=& \, L^{z}_{ijk} \, Q^{(4)}_{ijk},  
        \quad \Delta_{+} W^{(1)}_{ijk} &:=& \, L^{z}_{ijk} \, \left( {Q^{(1)}_{i j \, k+1}-Q^{(1)}_{ijk}} \right), 
        \quad  \Delta_{-} W^{(1)}_{ijk} &:=& \, L^{z}_{ijk} \, \left( {Q^{(1)}_{ijk}-Q^{(1)}_{i j \, k-1}} \right), \\
        W^{(2)}_{ijk} &:=& \, L^{z}_{ijk} \, Q^{(6)}_{ijk},  
        \quad \Delta_{+} W^{(2)}_{ijk} &:=& \, L^{z}_{ijk} \, \left( {Q^{(4)}_{i+1 \, j k}-Q^{(4)}_{ijk}} \right), 
        \quad  \Delta_{-} W^{(2)}_{ijk} &:=& \, L^{z}_{ijk} \, \left( {Q^{(4)}_{ijk}-Q^{(4)}_{i-1 \, j k}} \right),  \\
       W^{(3)}_{ijk} &:=& \, L^{z}_{ijk} \, Q^{(7)}_{ijk},  
        \quad \Delta_{+} W^{(3)}_{ijk} &:=& \, L^{z}_{ijk} \, \left( {Q^{(4)}_{i \, j+1 \, k}-Q^{(4)}_{ijk}} \right), 
        \quad  \Delta_{-} W^{(3)}_{ijk} &:=& \, L^{z}_{ijk} \, \left( {Q^{(4)}_{ijk}-Q^{(4)}_{i \, j-1 \, k}} \right), \\
        W^{(4)}_{ijk} &:=& \, L^{z}_{ijk} \, Q^{(10)}_{ijk},  
        \quad \Delta_{+} W^{(4)}_{ijk} &:=& \, L^{z}_{ijk} \, \left( {Q^{(4)}_{i j \, k+1}-Q^{(4)}_{ijk}} \right), 
        \quad  \Delta_{-} W^{(4)}_{ijk} &:=& \, L^{z}_{ijk} \, \left( {Q^{(4)}_{ijk}-Q^{(4)}_{i j \, k-1}} \right).      
\end{alignat}
\item Let $u^{(\ell)}_{ijk}$ be some component of $U^{(\ell)}_{ijk}$. Similarly,
	let $\Delta_{-} u^{(\ell)}_{ijk}$, $\Delta_{+} u^{(\ell)}_{ijk}$, $v^{(\ell)}_{ijk}$,
	$\Delta_{-} v^{(\ell)}_{ijk}$, $\Delta_{+} v^{(\ell)}_{ijk}$, $w^{(\ell)}_{ijk}$,
	$\Delta_{-} w^{(\ell)}_{ijk}$, and $\Delta_{+} w^{(\ell)}_{ijk}$ be 
	the corresponding components of $\Delta_{-} U^{(\ell)}_{ijk}$, $\Delta_{+} U^{(\ell)}_{ijk}$, $V^{(\ell)}_{ijk}$,
	$\Delta_{-} V^{(\ell)}_{ijk}$, $\Delta_{+} V^{(\ell)}_{ijk}$, $W^{(\ell)}_{ijk}$,
	$\Delta_{-} W^{(\ell)}_{ijk}$, and $\Delta_{+} W^{(\ell)}_{ijk}$.
\item Loop over all components and limit as follows:
\begin{algorithmic}
\STATE (a) Limit quadratic terms:
\begin{equation}
\begin{split}
      \tilde{u}^{(4)}_{ijk} := \text{mm} \left\{ u^{(4)}_{ijk}, \, \alpha_4 \,  \Delta_{-} u^{(4)}_{ijk}, \,
     \alpha_4 \,  \Delta_{+} u^{(4)}_{ijk}  \right\},  \quad
     \tilde{v}^{(4)}_{ijk} := \text{mm} \left\{ v^{(4)}_{ijk}, \, \alpha_4 \,  \Delta_{-} v^{(4)}_{ijk}, \,
     \alpha_4 \,  \Delta_{+} v^{(4)}_{ijk}  \right\},  \\ 
      \text{and} \quad
      \tilde{w}^{(4)}_{ijk} := \text{mm} \left\{ w^{(4)}_{ijk}, \, \alpha_4 \,  \Delta_{-} w^{(4)}_{ijk}, \,
     \alpha_4 \,  \Delta_{+} w^{(4)}_{ijk}  \right\}. \qquad \qquad \qquad \qquad
    \end{split}
\end{equation}
\STATE (b) If \, $\left|  \tilde{u}^{(4)}_{ijk} - u^{(4)}_{ijk} \right| + \left|  \tilde{v}^{(4)}_{ijk} - v^{(4)}_{ijk} \right| + \left| \tilde{w}^{(4)}_{ijk} - w^{(4)}_{ijk} \right| > 0$  \,
		or \, $\left|  {u}^{(4)}_{ijk} \right| = 0$ \, or \, $\left|  {v}^{(4)}_{ijk} \right| = 0$ \, or \, $\left|  {w}^{(4)}_{ijk} \right| = 0$ \,
	 then limit all the other terms ($\ell = 1,2,3$):
\begin{equation}
\begin{split}
      \tilde{u}^{(\ell)}_{ijk} := \text{mm} \left\{ u^{(\ell)}_{ijk}, \, \alpha_{\ell} \,  \Delta_{-} u^{(\ell)}_{ijk}, \,
     \alpha_{\ell} \,  \Delta_{+} u^{(\ell)}_{ijk}  \right\},  \quad
     \tilde{v}^{(\ell)}_{ijk} := \text{mm} \left\{ v^{(\ell)}_{ijk}, \, \alpha_{\ell} \,  \Delta_{-} v^{(\ell)}_{ijk}, \,
     \alpha_{\ell} \,  \Delta_{+} v^{(\ell)}_{ijk}  \right\},  \\ 
      \text{and} \quad
      \tilde{w}^{(\ell)}_{ijk} := \text{mm} \left\{ w^{(\ell)}_{ijk}, \, \alpha_{\ell} \,  \Delta_{-} w^{(\ell)}_{ijk}, \,
     \alpha_{\ell} \,  \Delta_{+} w^{(\ell)}_{ijk}  \right\}, \qquad \qquad \qquad \qquad
    \end{split}
\end{equation}
\STATE where in this work we take for $\alpha_{\ell}$ numbers in the range
	suggested by Krivodonova \cite{article:Kriv07}:  $\alpha_4 = \alpha_3 = \alpha_2 = \sqrt{\frac{3}{20}}$ \, and \, $\alpha_1 = \sqrt{\frac{1}{12}}$.
\end{algorithmic}
\item Convert characteristic values back to conservative values:
\begin{equation}
\begin{split}
     &{Q}^{(2)}_{ijk} \leftarrow R^{x}_{ijk} \, \tilde{U}^{(1)}_{ijk}, \quad 
     {Q}^{(3)}_{ijk} \leftarrow R^{y}_{ijk} \, \tilde{V}^{(1)}_{ijk}, \quad
     {Q}^{(4)}_{ijk} \leftarrow R^{z}_{ijk} \, \tilde{W}^{(1)}_{ijk}, \quad
     {Q}^{(5)}_{ijk} \leftarrow \text{mm} \left\{ R^{x}_{ijk} \, \tilde{U}^{(2)}_{ijk}, \, R^{y}_{ijk} \, \tilde{V}^{(2)}_{ijk} \right\}, \\
     &{Q}^{(6)}_{ijk} \leftarrow \text{mm} \left\{ R^{x}_{ijk} \, \tilde{U}^{(3)}_{ijk}, \, R^{z}_{ijk} \, \tilde{W}^{(2)}_{ijk} \right\}, \quad
     {Q}^{(7)}_{ijk} \leftarrow \text{mm} \left\{ R^{y}_{ijk} \, \tilde{V}^{(3)}_{ijk}, \, R^{z}_{ijk} \, \tilde{W}^{(3)}_{ijk} \right\}, \quad
     {Q}^{(8)}_{ijk} \leftarrow R^{x}_{ijk} \, \tilde{U}^{(4)}_{ijk}, \\ 
     &{Q}^{(9)}_{ijk} \leftarrow R^{y}_{ijk} \, \tilde{V}^{(4)}_{ijk}, \quad
     {Q}^{(10)}_{ijk} \leftarrow R^{z}_{ijk} \, \tilde{W}^{(4)}_{ijk}.
      \end{split}
\end{equation}
\end{enumerate}

\section{Semi-discrete entropy stability theory for DG applied to ideal MHD}
\label{sec:dg_stability}
As shown in Section \ref{sec:symmhyp}, the ideal MHD equations
cannot be put in symmetric hyperbolic form without adding to the
equations a term that is proportional to the divergence of the magnetic
field. As it turns out, it is precisely this fact that causes many standard
high-resolution methods such as wave propagation, weighted
essentially non-oscillatory, and DG finite element methods to be
numerically unstable if the divergence of the magnetic field is
not properly controlled. In this section we briefly review the
key theorem from Barth \cite{article:Barth05} that rigorously
proves this statement. 
 
For simplicity of exposition we will not give consider the entropy
stability for the DG method represented by 
\eqref{eqn:semidiscrete_fem1}--\eqref{eqn:semidiscrete_fem2},
but instead consider a modification of this method where it is
the entropy variables \eqref{eqn:entropy_vars} that are expanded as a sum of 
polynomials rather than the conserved variables:
\begin{equation}
  v^h\left(t, x, y, z \right) \Bigl|_{\Tm_i} :=
  \sum_{m=1}^{M(M+1)(M+2)/6} V^{(m)}_i(t) \, \varphi^{(m)}\left(\xi, \eta, \zeta \right).
\end{equation}
We emphasize that this modification is only for the discussion in this section,
and that the constrained transport method we describe in Section \ref{sec:ct_scheme}
is based on using conserved variables via \eqref{eqn:semidiscrete_fem1}--\eqref{eqn:semidiscrete_fem2}.
Using $v^h$ instead of $q^h$, we arrive at the following
semi-discrete variational problem for all $t\in[0,T]$:
\begin{equation}
\begin{split}
     \text{Find} \quad &v^h(t,{\bf x}) \in \WS^{h,8}_{\deg} \quad \text{such that} \\
     &\tilde{B}_{\text{DG}}\left(v^h,w^h \right) = 0 \qquad
       \text{for all} \quad w^h(t,{\bf x}) \in \WS^{h,8}_{\deg},
       \end{split}
\label{eqn:semidiscrete_entropy1}
\end{equation}
where
\begin{equation}
\begin{split}
   \widetilde{B}_{\text{DG}}\left(v^h,w^h \right) :=  &\sum_{i=1}^N \ \iiint_{\Tm_i}  \left\{ w^h \cdot
q(v^h)_{,t} - \nabla w^h : {\bf F}(v^h) + \sigma \, w^h \cdot \chi_{,v} \, \nabla \cdot \B \ \right\} \, d{\bf x} \\
 +   &\sum_{i=1}^N  \varoiint_{\partial \Tm_i}
       w^h_{-}  \cdot \NF\left(v^h_{-}({\bf s}), \, v^h_{+}({\bf s}); \, {\bf n}\right) \, 
       d{\bf s}.
\label{eqn:semidiscrete_entropy2}
\end{split}
\end{equation}
We view each of the conservative variables, $\rho$, $\rho \u$, ${\mathcal E}$, and $\B$, as functions of $v^h$.
Motivated by the entropy analysis from Section \ref{sec:symmhyp},  the above expression includes a term proportional to the divergence of the magnetic field with a tunable parameter $\sigma$: if $\sigma\equiv0$ the additonal term vanishes,
if $\sigma\equiv1$ then we have the full Powell source term from \eqref{eqn:sym_hyp}.

\begin{thm}[
Barth \cite{article:Barth05}]
\label{thm:entropy_stability}
Let $v^h(t,{\bf x})$ be the solution to the semi-discrete DG variational problem 
\eqref{eqn:semidiscrete_entropy1}--\eqref{eqn:semidiscrete_entropy2}.
Assume that all of the following conditions are satisfied:
\begin{enumerate}
\item \underline{Powell term or div-free on element:} \quad   $\sigma \equiv 1$ \quad or  \quad $\left(\nabla \cdot \B \right)\bigl|_{\Tm_i} = 0$ \quad on each element $\Tm_i$;
\item \underline{Continuous normal component across edges:} \quad ${\bf n}_e \cdot \left( {\bf B}_{+} -  {\bf B}_{-} \right) = 0$
\quad across each face $e$ (or edge in 2D).
\item \underline{System E-flux condition:}
\begin{equation}
\label{eqn:system_eflux}
\llbracket v^h \rrbracket^{+}_{-} \cdot 
\left\{ {\mathcal F}\left(v^h_{-}, \, v^h_{+}; \, {\bf n} \right) - 
{\bf F}\left(v^h_{-} + \theta \, \llbracket v^h \rrbracket^{+}_{-} \right) \cdot {\bf n} \right\} \le 0, \quad \forall \, \theta \in [0,1],
\end{equation}
where $\llbracket v^h \rrbracket^{+}_{-} := v^h_{+} - v^h_{-}$.
\end{enumerate}
The numerical solution $v^h(t,{\bf x})$ then satisfies the following entropy inequality on
each element $\Tm_i$:
\begin{equation}
\label{eqn:semidiscrete_entropy_condition}
   \frac{d}{d t} \iiint_{\Tm_i}  U(v^h) \, d{\bf x} + \varoiint_{\partial \Tm_i} {\mathcal G}(v^h_{-}, \, v^h_{+}; \,
   {\bf n}) \, ds \le 0,
\end{equation}
where $U$ is the entropy function and ${\mathcal G}(v^h_{-}, \, v^h_{+}; \,
   {\bf n})$ is the numerical entropy flux:
\begin{equation}
\label{eqn:num_entropy_flux}
   {\mathcal G}\left(v_{-}^h, \, v_{+}^h, \, {\bf n}\right):=
 \big\langle v^h \big\rangle_{-}^{+} \cdot
    {\mathcal F}\left(v_{-}^h, \, v_{+}^h; \, {\bf n}\right) - 
  \big\langle  \dualeflux( v^h ) \cdot {\bf n} - \chi( v^h ) \, \B( v^h ) \cdot {\bf n}
     \big\rangle_{-}^{+}, 
\end{equation}
where $\big\langle \cdot \big\rangle^+_{-}$ denotes the average of the values
on either side of the element face (or edge in 2D) and $\dualeflux$ is the dual entropy flux \eqref{eqn:dual_entropy_flux}.
This local entropy inequality also implies a global semi-discrete entropy inequality:
\begin{equation}
\label{eqn:semidiscrete_entropy_condition_global}
   \frac{d}{d t} \iiint_{\Omega}  U\left(v^h \right) \, d{\bf x} \le 0.
\end{equation}

\end{thm}

The proof of this theorem is given in Barth \cite{article:Barth05} and follows from 
evaluating the bilinear form in \eqref{eqn:semidiscrete_entropy1}--\eqref{eqn:semidiscrete_entropy2}
with test functions set to the entropy variables:  $\tilde{B}_{\text{DG}}\left( v^h, v^h \right) = 0$. 
This is possible because in the formulation of \eqref{eqn:semidiscrete_entropy1}--\eqref{eqn:semidiscrete_entropy2},
$v^h$ and $w^h$ are both in the same finite element space: $\WS^{h,8}_{\deg}$.
After application of various identities involving entropy variables and dual entropy variables,
as well as the assumptions on the magnetic field, the semi-discrete entropy inequalities
\eqref{eqn:semidiscrete_entropy_condition}  and \eqref{eqn:semidiscrete_entropy_condition_global} follow.

\begin{remark}
The numerical method proposed in the next section (\S \ref{sec:ct_scheme}) differs from the
DG method used in Theorem \ref{thm:entropy_stability} in several important ways: (1)
the conserved variables are used as trial variables, not the entropy variables; (2)
a strong-stability-preserving Runge-Kutta (SSP-RK) method is used to turn
the semi-discrete equations into fully discrete equations; and (3) the volume and surface
integrals are replaced by appropriate Gaussian quadrature rules.
Despite these differences, the core result of Theorem \ref{thm:entropy_stability}
extends to SSP-RK DG schemes: a globally divergence-free magnetic field
leads to entropy-stabilized discretizations. As a matter of practice we
rely on the limiters (see \S \ref{sec:limiters}) to inject additional numerical dissipation
in order to overcome the entropy errors introduced through the use of conserved
variables, SSP-RK time-stepping, and
Gaussian quadrature rules.
\end{remark}

In the next section (\S \ref{sec:ct_scheme})
we develop a constrained transport SSP-RK DG scheme. In Section \ref{sec:numex}
we demonstrate through numerical examples that this approach
does not adversely affect the convergence rate of the overall scheme and
that, in conjunction with moment-based limiters, produces numerically
stable solutions of ideal MHD even in the presence of shocks.

\section{Globally divergence-free constrained transport for DG-FEMs}
\label{sec:ct_scheme}
The lesson that should be taken from Theorem \ref{thm:entropy_stability}
is that we need to control two quantities simultaneously:
(1) the divergence of the magnetic field within each element, as well
as (2) the jump in the normal components across each element edge.
This, of course, means that we must control the {\it global divergence}
of the magnetic field on the computational grid. 
In this section we provide such a strategy both for 2D and 3D Cartesian meshes. 

\subsection{Basic constrained transport algorithm}
In order to compute a globally divergence-free magnetic field at each time step,
we introduce the magnetic vector potential,
\begin{equation}
   \B = \nabla \times \A,
\end{equation}
which satisfies the magnetic induction equation:
\begin{equation}
\label{eqn:induction}
  \A_{,t}  = \u \times \B,
\end{equation}
where we have assumed the Weyl gauge condition 
(see Helzel et al. \cite{article:HeRoTa11} for a discussion
of various gauge conditions).

We summarize the proposed constrained transport below by presenting a 
single forward Euler time-step. 
Extension to high-order strong stability-preserving Runge-Kutta (SSP-RK) methods is straightforward,
since SSP-RK time-stepping methods are  convex combinations of forward Euler steps. A 
single time-step of the proposed CT method from time $t=t^n$ to time $t=t^{n+1}$ consists of the following sub-steps:

\medskip

\begin{description}
  \item[{\bf Step 0.}] The following state variables are given as initial data at time $t=t^n$:
  \begin{align}
   \text{Mass, momentum, energy:} & \quad
   		   \rho^h \left( t^n, {\bf x} \right) \in \WS^{h,1}_{\deg}, 
		\quad \rho \u^h \left( t^n, {\bf x} \right) \in \WS^{h,3}_{\deg},
		\quad \En^h \left( t^n, {\bf x} \right) \in \WS^{h,1}_{\deg}, \\
   \text{Magnetic field:} & \quad
         \B^h \left( t^n, {\bf x} \right) \in \WS^{h,3}_{{\deg}+1}, \quad \widetilde{\B}^h \left( t^n, {\bf x} \right) \in \WS^{h,3}_{{\deg}}, \\
   \text{Magnetic potential:} & \quad
        \A^h \left( t^n, {\bf x} \right) \in \WS^{h,3}_{{\deg}}.
  \end{align}
   Note that there are two representations of the magnetic field.
               One is globally divergence-free and lives in a broken finite element space
               of one degree higher than all the other conserved variables: $\B^h\left( t^n, {\bf x} \right) \in \WS^{h,3}_{{\deg}+1}$.
               The other is not globally divergence-free and lives in the same
               broken finite element space as all the other conserved variables: $\widetilde{\B}^h\left( t^n, {\bf x} \right) \in \WS^{h,3}_{{\deg}}$.
               
\medskip
               	
  \item[{\bf Step 1.}] Update $\rho^h$, $\rho \u^h$, $\En^h$, $\widetilde{\B}^h$ using the DG scheme
  with a forward Euler time-step:
	           \begin{equation}
	           	\frac{1}{|\Tm_i|} \iiint_{\Tm_i} w^h_i \cdot \begin{bmatrix}
			     \rho^h\left(t^{n+1}, {\bf x}\right) \\
			     \rho \u^h\left(t^{n+1}, {\bf x}\right) \\
			     \En^h\left(t^{n+1}, {\bf x}\right) \\
			     \widetilde{\B}^h\left(t^{n+1}, {\bf x}\right)
			\end{bmatrix} \, d\xv= 
			\frac{1}{|\Tm_i|} \iiint_{\Tm_i} w^h_i \cdot \begin{bmatrix}
			     \rho^h\left(t^n, {\bf x}\right) \\
			     \rho \u^h\left(t^n, {\bf x}\right) \\
			     \En^h\left(t^n, {\bf x}\right) \\
			     \widetilde{\B}^h\left(t^n, {\bf x}\right)
			\end{bmatrix} \, d\xv
			  + \Delta t \, {\mathcal L}_i \left(
			  \begin{bmatrix}
			     \rho^h\left(t^n, {\bf x}\right) \\
			     \rho \u^h\left(t^n, {\bf x}\right) \\
			     \En^h\left(t^n, {\bf x}\right) \\
			     \B^h\left(t^n, {\bf x}\right)
			\end{bmatrix}, w^h_i \right).
	           \end{equation}
	The updated magnetic field, $\widetilde{\B}^h(t^{n+1}, {\bf x})$, is {\it not} globally divergence-free;
	we view this value as the {\it predicted} magnetic field.
	Note that the operator ${\mathcal L}_i$ is evaluated using the globally divergence-free
	magnetic field: ${\B}^h(t^{n}, {\bf x})$; this is important for numerical stability.

\medskip

 \item[{\bf Step 2.}] Update the  magnetic potential using a forward Euler time-step on the induction equation:
	\begin{equation}
	\hspace{-5mm}
		 \frac{1}{|\Tm_i|} \iiint_{\Tm_i} w^h_i \cdot  \A^h\left(t^{n+1}, {\bf x}\right) \, d\xv = 
		\frac{1}{|\Tm_i|} \iiint_{\Tm_i} w^h_i \cdot  \A^h\left(t^n, {\bf x}\right) \, d\xv  + \frac{{\Delta t}}{|\Tm_i|}
		 \, \iiint_{\Tm_i} w^h_i \cdot
		\left( \frac{\rho \u^h\left(t^n, {\bf x} \right) \times \widetilde{\B}^h(t^{n}, {\bf x})}{\rho^h\left(t^{n}, \xv\right)}  \right) \, d\xv.
	\end{equation}
	
	\medskip
	
   \item[{\bf Step 3.}] From the magnetic potential and the predicted magnetic field,
   	construct a globally divergence-free magnetic field:
	\begin{equation}
	\label{eqn:divfree_construct}
		 \B^h(t^{n+1}, {\bf x}) = \text{Div-Free-Construct}\left( \widetilde{\B}^h(t^{n+1}, {\bf x}), \,
		 	 \A^h(t^{n+1}, {\bf x}) \right).
	\end{equation}
	 The details of the `Div-Free-Construct' operator are described for the 2D case in \S \ref{sec:curl2d}
	and the 3D case in \S \ref{sec:curl3d}.
	
	\medskip
	
  \item[{\bf Step 4.}] Synchronize the two versions of the magnetic field by performing a simple $L^2$-projection
  from $\WS^{h,3}_{\deg+1}$ to $\WS^{h,3}_{\deg}$:
  \begin{equation}
  	\widetilde{\B}^h \left( t^{n+1}, {\bf x} \right) = \text{$L^{2}$--Project}\left( \B^h(t^{n+1}, {\bf x}) \right).
  \end{equation}
  Note that in the case of a multi-stage Runge-Kutta scheme, we only perform this synchronization
  step once per time-step (i.e., not after every stage); that is, we only perform the synchronization at the end of the full 
  time-step, once all the stages have been completed. We have found that synchronizing only
  once per time step significantly reduces unphysical oscillations in the case of
   solutions with shocks.
\end{description}

\begin{figure}[!t]
\begin{center}
\begin{tabular}{cc}
  (a) \includegraphics[height=45mm]{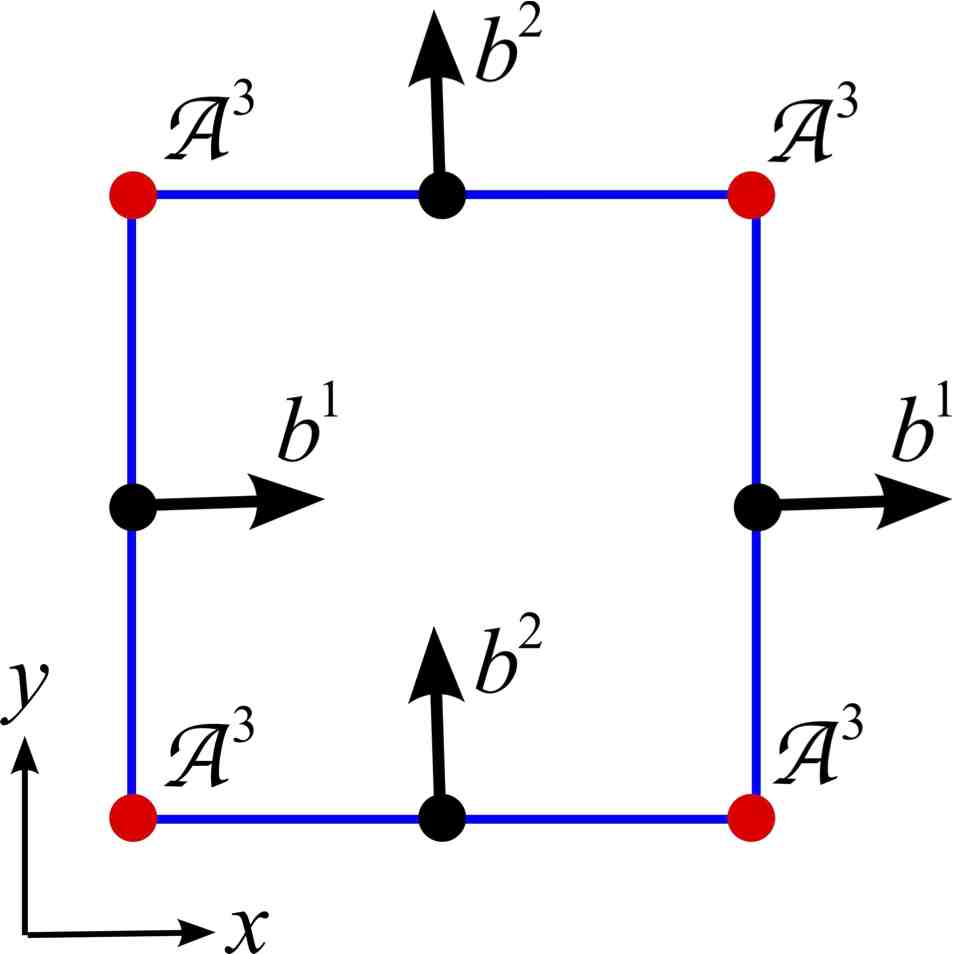} &
  (b) \includegraphics[height=45mm]{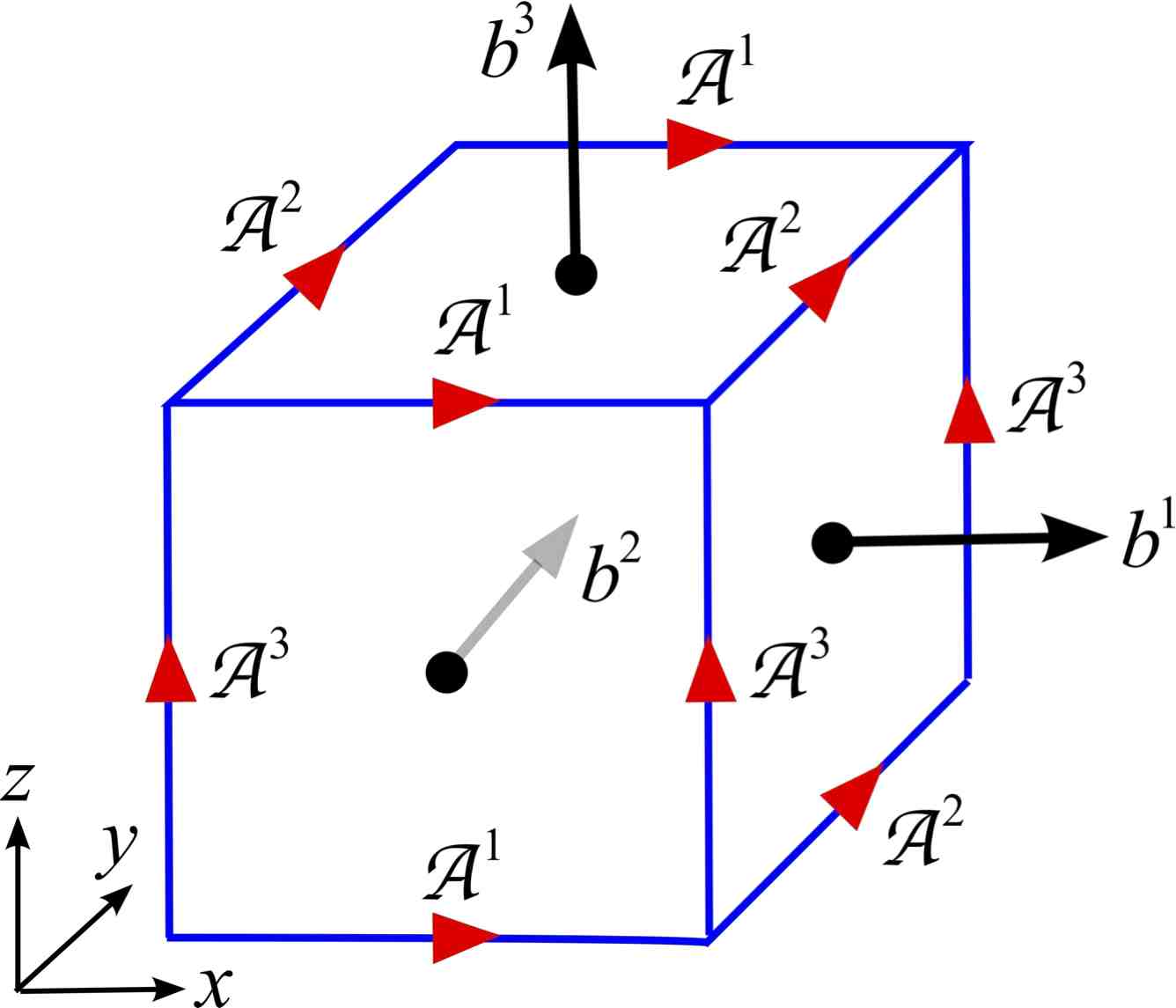} 
\end{tabular}
  \caption{Staggered magnetic field and magnetic potential configurations for Cartesian mesh elements.
  Shown in Panel (a) is a 2D mesh element with the magnetic field components on mesh
  edges and the potential on mesh corners. Shown in Panel (b) is a 3D mesh element with the magnetic 
  field components on mesh
  faces and the potential on mesh edges. 
  \label{fig:yee}}
\end{center}
\end{figure}

\subsection{2D construction of a divergence-free magnetic field}
\label{sec:curl2d}
The final pieces missing from the constrained transport method
as proposed in this work are the details of the divergence-free
construction step as referred to in \eqref{eqn:divfree_construct}.
The approach advocated here is a modification of similar 
divergence-free constructions used by Balsara \cite{article:Ba04}
and Li, Xu, and Yakovlev \cite{article:LiXuYa11}.
In fact, the divergence-free reconstruction of the magnetic
field is directly related to ideas from discrete exterior calculus and, in particular,
to Whitney forms (e.g., see Bossavit \cite{article:Bossavit88}). Despite these connections, we do need
the full machinery of discrete exterior calculus and Whitney forms in order to
develop the proposed numerical scheme.

In 2D the divergence-free reconstruction makes use of 
magnetic field components normal to element edges and magnetic vector
potential values on element corners. This grid staggering 
is illustrated in Figure \ref{fig:yee}(a).
We give the full details of the divergence-free reconstruction below:
\begin{description}
  \item[{\bf Step 0.}] Start with the {\it predicted} magnetic field and the magnetic potential:
  \begin{equation}
	 \widetilde{\B}^h \left( {\bf x} \right) \cdot {\bf e}^m  \Bigl|_{\Tm_{ij}}
	 = \sum_{\ell=1}^{6}  \Bt^{ \, m(\ell)}_{ij} \, \varphi^{(\ell)}\left( \xi, \eta \right) \qquad
	 \text{and} \qquad
  	\A^h \left( {\bf x} \right) \cdot {\bf e}^3  \Bigl|_{\Tm_{ij}} 
		= \sum_{\ell=1}^{6}  A^{3(\ell)}_{ij} \, \varphi^{(\ell)}\left( \xi, \eta \right).
  \end{equation}
   \item[{\bf Step 1.}] Interpolate the magnetic potential to mesh corners:
\begin{equation}
{\mathcal A}^3_{i-\half \, j-\half} := \frac{1}{4} \sum_{\ell=1}^{6} \biggl[  
		A^{3(\ell)}_{ij} \, \varphi^{(\ell)}{\left( -1, -1 \right)} + 
		A^{3(\ell)}_{ij-1} \, \varphi^{(\ell)}\left( -1,  1 \right)   + 
		A^{3(\ell)}_{i-1j} \, \varphi^{(\ell)}\left(  1,  -1 \right) +
		A^{3(\ell)}_{i-1 j-1} \, \varphi^{(\ell)}\left(  1, 1 \right)
	 \biggr].
\end{equation}
 \item[{\bf Step 2.}] On each element edge, define a DG representation for the normal components
		of the magnetic field:
		\begin{equation}
		\label{eqn:edge_2d}
		b^1_{i-\half \, j} = \sum_{\ell=1}^{3}  b^{1(\ell)}_{i-\half \, j} \, \varphi_{\text{1D}}^{(\ell)}\left( \alpha  \right)
		\qquad \text{and} \qquad
		b^2_{i \, j-\half} = \sum_{\ell=1}^{3}  b^{2(\ell)}_{i \, j-\half} \, \varphi_{\text{1D}}^{(\ell)}\left( 
		\alpha \right),
		\end{equation}
		where $\varphi_{\text{1D}}^{(\ell)}$
		 are the 1D Legendre polynomials.
		The coefficients of $b^1$ and $b^2$ are primarily computed from direct interpolation of
		element centered magnetic field values: $\widetilde{\B}^h$. The exceptions
		to this are the average magnetic values on the edges, which, in order to guarantee
		zero divergence, are computed from finite differences of the magnetic potential
		on element corners. The detailed equations can be written as follows:
		\begin{align}
		\label{eqn:ave_b_2d}
 b^{1(1)}_{i-\half \, j} &= \frac{{\mathcal A}^3_{\, i-\half \, j+\half} -   {\mathcal A}^3_{\, i-\half \, j-\half}}{\Delta y},  & 
 b^{2(1)}_{i \, j-\half} &= \frac{{\mathcal A}^3_{\, i-\half \, j-\half} -   {\mathcal A}^3_{\, i+\half \, j-\half}}{\Delta x}, \\
 b^{1(2)}_{i-\half \, j} &=  \frac{\Bt^{ \, 1(3)}_{i-1 j} + \sqrt{3} \, \Bt^{ \, 1(4)}_{i-1 j} +
 	\Bt^{ \, 1(3)}_{ij} - \sqrt{3} \, \Bt^{ \, 1(4)}_{ij}}{2},  &
 b^{2(2)}_{i \, j-\half} &=  \frac{\Bt^{ \, 2(2)}_{ij-1} + \sqrt{3} \, \Bt^{ \, 2(4)}_{ij-1} +
 	\Bt^{2(2)}_{ij} - \sqrt{3} \, \Bt^{ \, 2(4)}_{ij}}{2}, \\
 b^{1(3)}_{i-\half \, j} &= \frac{\Bt^{ \, 1(6)}_{i-1 j} + \Bt^{ \, 1(6)}_{ij}}{2},  &
 b^{2(3)}_{i \, j-\half} &= \frac{\Bt^{ \, 2(5)}_{i j-1} + \Bt^{ \, 2(5)}_{ij}}{2}.
\end{align}
 \item[{\bf Step 3.}] The final step is to construct globally divergence-free magnetic field values:
\begin{equation}
\label{eqn:bfield_full_2d}
{\B}^h \left( {\bf x} \right) \cdot {\bf e}^1  \Bigl|_{\Tm_{ij}}
	 = \sum_{\ell=1}^{10}  B^{ \, 1(\ell)}_{ij} \, \varphi^{(\ell)}\left( \xi, \eta \right) \qquad
	 \text{and} \qquad
	 {\B}^h \left( {\bf x} \right) \cdot {\bf e}^2  \Bigl|_{\Tm_{ij}}
	 = \sum_{\ell=1}^{10}  B^{ \, 2(\ell)}_{ij} \, \varphi^{(\ell)}\left( \xi, \eta \right).
	 \end{equation}
This is achieved by enforcing the following conditions: 
\begin{enumerate}
\item ${\B}^h \left( {\bf x} \right) \cdot {\bf e}^1$ exactly
matches $b^1$ from \eqref{eqn:edge_2d} on the left $(i-\half,j)$ and right $(i+\half,j)$ edges;
\item ${\B}^h \left( {\bf x} \right) \cdot {\bf e}^2$ exactly
matches $b^2$ from \eqref{eqn:edge_2d} on the bottom $(i, j-\half)$ and top $(i,j+\half)$ edges;
\item $\nabla \cdot {\B}^h \Bigl|_{\Tm_{ij}} = \text{constant}$; and
\item Any coefficients in \eqref{eqn:bfield_full_2d} that remain as free parameters are set to zero.
\end{enumerate}
The detailed equations for the 1-component can be written as follows:
\begin{equation}
\label{eqn:div-free-b-2d-1}
\begin{split}
B^{1(1)}_{ij} &= \half \left(\ar{1}+\al{1}\right) + \halfsqt \dxody \left(\br{2}-\bl{2} \right), \\
B^{1(2)}_{ij} &= \halfsqt \left(\ar{1}-\al{1}\right)+\halfsqft \dxody \left(\br{3}-\bl{3}\right), \\
B^{1(3)}_{ij} &= \half \left(\ar{2}+\al{2}\right), \quad
B^{1(4)}_{ij} = \halfsqt \left(\ar{2}-\al{2}\right), \\
B^{1(5)}_{ij} &= \halfsqft \dxody \left(\bl{2}-\br{2}\right),  \quad
B^{1(6)}_{ij} =\half \left(\ar{3}+\al{3}\right), \quad
B^{1(7)}_{ij} = 0, \\
B^{1(8)}_{ij} &= \halfsqt \left(\ar{3}-\al{3}\right), \quad
B^{1(9)}_{ij} = \halfsqtf \dxody \left(\bl{3}-\br{3}\right), \quad
B^{1(10)}_{ij} = 0,
\end{split}
\end{equation}
The detailed equations for the 2-component can be written as follows:
\begin{equation}
\label{eqn:div-free-b-2d-2}
\begin{split}
B^{2(1)}_{ij} &= \half \left(\br{1}+\bl{1}\right)+\halfsqt \dyodx \left(\ar{2}-\al{2}\right), \quad
B^{2(2)}_{ij} = \half \left(\br{2}+\bl{2}\right), \\
B^{2(3)}_{ij} &= \halfsqt \left(\br{1}-\bl{1}\right)+\halfsqft \dyodx \left(\ar{3}-\al{3}\right), \quad
B^{2(4)}_{ij} = \halfsqt \left(\br{2}-\bl{2}\right), \\ 
B^{2(5)}_{ij} &= \half \left(\br{3}+\bl{3}\right), \quad
B^{2(6)}_{ij} =  \halfsqft \dyodx \left(\al{2}-\ar{2}\right), \\
B^{2(7)}_{ij} &=  \halfsqt \left(\br{3}-\bl{3}\right), \quad
B^{2(8)}_{ij} =  B^{2(9)}_{ij} =  0, \quad
B^{2(10)}_{ij} =  \halfsqtf \dyodx \left(\al{3}-\ar{3}\right).
\end{split}
\end{equation}
\end{description}

A straightforward calculation reveals that the (pointwise) divergence
of a magnetic field of the form \eqref{eqn:bfield_full_2d} with coefficients 
given by \eqref{eqn:div-free-b-2d-1} and \eqref{eqn:div-free-b-2d-2}  is
\begin{equation}
\label{eqn:div_zero_2d}
 \nabla \cdot \B^h \Bigl|_{\Tm_{ij}} = 
 \frac{2}{\Delta x} \frac{\partial}{\partial \xi} \left( {\B}^h \left( {\bf x} \right) \cdot {\bf e}^1  \Bigl|_{\Tm_{ij}} 
\right)
+ \frac{2}{\Delta y} \frac{\partial}{\partial \eta} 
\left( {\B}^h \left( {\bf x} \right) \cdot {\bf e}^2  \Bigl|_{\Tm_{ij}} 
\right) =
 \frac{\ar{1}-\al{1}}{\Delta x} +\frac{\br{1}-\bl{1}}{\Delta y},
\end{equation}
which is constant. Using definitions
\eqref{eqn:ave_b_2d}, this divergence can be shown to vanish:
\begin{equation*}
\nabla \cdot \B^h \Bigl|_{\Tm_{ij}} =   \frac{ 
  {\mathcal A}^3_{\, i+\half \, j+\half} - {\mathcal A}^3_{\, i+\half \, j-\half}
- {\mathcal A}^3_{\, i-\half \, j+\half} + {\mathcal A}^3_{\, i-\half \, j-\half} 
+ {\mathcal A}^3_{\, i-\half \, j+\half} - {\mathcal A}^3_{\, i+\half \, j+\half}
- {\mathcal A}^3_{\, i-\half \, j-\half} + {\mathcal A}^3_{\, i+\half \, j-\half} 
}{\Delta x \, \Delta y} = 0.
\end{equation*}
We note that we have actually achieved a globally divergence free
magnetic field, $\B^h(\xv)$, since we have the two sufficient ingredients:
\begin{enumerate}
\item $\B^h(\xv)$ restricted to the interior of each element is exactly
divergence-free (see \eqref{eqn:div_zero_2d}); and
\item The normal components of $\B^h(\xv)$ are continuous across
each element edge (see \eqref{eqn:div-free-b-2d-1}--\eqref{eqn:div-free-b-2d-2}).
\end{enumerate}

\subsection{3D construction of a divergence-free magnetic field}
\label{sec:curl3d}
The same basic principle used in 2D can be extended to construct a globally
divergence-free magnetic field in 3D. 
In 3D the divergence-free reconstruction makes use of 
magnetic field components normal to element faces and magnetic vector
potential values on element edges. This grid staggering 
is illustrated in Figure \ref{fig:yee}(b).
We outline the divergence-free reconstruction procedure 
below:

\smallskip

\begin{description}
  \item[{\bf Step 0.}] Start with the {\it predicted} magnetic field and magnetic potential:
  \begin{equation}
	 \widetilde{\B}^h \left( {\bf x} \right) \cdot {\bf e}^m  \Bigl|_{\Tm_{ijk}}
	 = \sum_{\ell=1}^{10}  \Bt^{ \, m(\ell)}_{ijk} \, \varphi^{(\ell)}\left( \xi, \eta, \zeta \right)
	 \qquad \text{and} \qquad
  	\A^h \left( {\bf x} \right) \cdot {\bf e}^m  \Bigl|_{\Tm_{ijk}} 
		= \sum_{\ell=1}^{10}  A^{m(\ell)}_{ijk} \, \varphi^{(\ell)}\left( \xi, \eta, \zeta \right).
  \end{equation}
  
  \medskip
  
   \item[{\bf Step 1.}] Interpolate the magnetic potential to mesh edges:
\begin{gather}
\begin{split}
{\mathcal A}^1_{i \, j-\half \, k-\half} := \frac{1}{8} \sum_{\ell=1}^{10} \int_{-1}^{1} \biggl[  
		&A^{1(\ell)}_{ijk} \, \varphi^{(\ell)}\left( \xi, -1, -1 \right) + 
		A^{1(\ell)}_{ij k-1} \, \varphi^{(\ell)}\left( \xi, -1,  1  \right) \\ + \, \, 
		& A^{1(\ell)}_{i j-1 k} \, \varphi^{(\ell)}\left( \xi, 1,  -1 \right) +
		A^{1(\ell)}_{i j-1 k-1} \, \varphi^{(\ell)}\left(  \xi, 1, 1 \right)
	 \biggr] \, d\xi,
\end{split} \\
\begin{split}
{\mathcal A}^2_{i-\half \, j \, k-\half} := \frac{1}{8}  \sum_{\ell=1}^{10} \int_{-1}^{1} \biggl[  
		&A^{2(\ell)}_{ijk} \, \varphi^{(\ell)}\left( -1, \eta, -1 \right) + 
		A^{2(\ell)}_{i j k-1} \, \varphi^{(\ell)}\left( -1,  \eta, 1 \right) \\ + \, \, 
		& A^{2(\ell)}_{i-1 j k} \, \varphi^{(\ell)}\left(  1,  \eta, -1 \right) +
		A^{2(\ell)}_{i-1 j k-1} \, \varphi^{(\ell)}\left(  1, \eta, 1 \right)
	 \biggr] \, d\eta,
\end{split} \\
\begin{split}
{\mathcal A}^3_{i-\half \, j-\half \, k} := \frac{1}{8}   \sum_{\ell=1}^{10}  \int_{-1}^{1} \biggl[  
		&A^{3(\ell)}_{ijk} \, \varphi^{(\ell)}\left( -1, -1, \zeta \right) + 
		A^{3(\ell)}_{ij-1 k} \, \varphi^{(\ell)}\left( -1,  1, \zeta \right) \\ + \, \, 
		& A^{3(\ell)}_{i-1j k} \, \varphi^{(\ell)}\left(  1,  -1, \zeta \right) +
		A^{3(\ell)}_{i-1 j-1 k} \, \varphi^{(\ell)}\left(  1, 1, \zeta \right)
	 \biggr] \, d\zeta.
\end{split}
\end{gather}
Note that these magnetic potential values are the edge averages of the magnetic potential;
for example, ${\mathcal A}^1_{i \, j-\half \, k-\half}$ is an approximation to the average of 
$A^1(x,y_j - \Delta y/2, z_k - \Delta z/2)$ over the interval $x \in [x_i - \Delta x/2, x_i + \Delta x/2]$.
 
  \medskip
  
\item[{\bf Step 2.}] On each element face, define a DG representation for the normal components
		of the magnetic field:
		\begin{align}
		\label{eqn:face_3d}
		\left\{ b^1_{i-\half \, j \, k}, \, b^2_{i \, j-\half \, k}, \,  
		b^3_{i \, j \, k-\half} \right\} &= \sum_{\ell=1}^{6}  
		\left\{ b^{1(\ell)}_{i-\half \, j \, k}, \, b^{2(\ell)}_{i \, j-\half \, k}, \,
		  b^{3(\ell)}_{i \, j \, k-\half} \right\}
			\, \varphi_{\text{2D}}^{(\ell)}\left( \alpha, \, \beta \right),
		\end{align}
		where
		$\varphi_{\text{2D}}^{(\ell)}$
	        are the 2D Legendre polynomials \eqref{eqn:test-fun-cart}.
		The coefficients of $b^1$, $b^2$, and $b^3$ are primarily computed from direct interpolation of
		element centered magnetic field values, $\widetilde{\B}^h$. The exceptions
		to this are the average magnetic values on the faces, which, in order to guarantee
		zero divergence, are computed from finite differences of the average magnetic potential
		on element edges. The full formulas are given by \eqref{eqn:ave_b_3d_1}-\eqref{eqn:ave_b_3d_3}.

  \medskip
  
 \item[{\bf Step 3.}] The final step is to construct globally divergence-free magnetic field values:
\begin{equation}
\label{eqn:bfield_full_3d}
{\B}^h \left( {\bf x} \right) \cdot {\bf e}^m  \Bigl|_{\Tm_{ijk}}
	 = \sum_{\ell=1}^{20}  B^{ \, m(\ell)}_{ijk} \, \varphi^{(\ell)}\left( \xi, \eta, \zeta \right).
	 \end{equation}
This is achieved by enforcing the following conditions: 
\begin{enumerate}
\item ${\B}^h \left( {\bf x} \right) \cdot {\bf e}^1$ exactly
matches $b^1$ from \eqref{eqn:face_3d} on the left $(i-\half,j,k)$ and right $(i+\half,j,k)$ faces;
\item ${\B}^h \left( {\bf x} \right) \cdot {\bf e}^2$ exactly
matches $b^2$ from \eqref{eqn:face_3d} on the front $(i, j-\half,k)$ and back $(i,j+\half,k)$ faces;
\item ${\B}^h \left( {\bf x} \right) \cdot {\bf e}^3$ exactly
matches $b^3$ from \eqref{eqn:face_3d} on the bottom $(i,j,k-\half)$ and top $(i,j,k+\half)$ faces;
\item $\nabla \cdot {\B}^h \Bigl|_{\Tm_{ijk}} = \text{constant}$; and
\item Any coefficients in \eqref{eqn:bfield_full_3d} that remain as free parameters are set to zero.
\end{enumerate}
The full formulas are given by \eqref{eqn:div-free-b-3d-1}--\eqref{eqn:div-free-b-3d-3}.
\end{description}

\medskip

A straightforward calculation reveals that the (pointwise) divergence
of a magnetic field of the form \eqref{eqn:bfield_full_3d} with coefficients 
given by \eqref{eqn:div-free-b-3d-1}--\eqref{eqn:div-free-b-3d-3-end} is
\begin{equation}
 \nabla \cdot \B^h \Bigl|_{\Tm_{ijk}} =
 \frac{b^1_{i+\half j k}-b^1_{i-\half j k}}{\Delta x} +\frac{b^2_{i j+\half k}-b^2_{i j-\half k}}{\Delta y}
 +\frac{b^3_{i j k+\half}-b^3_{i j k-\half}}{\Delta z},
\end{equation}
which is constant. Using definitions
\eqref{eqn:ave_b_3d_1}, \eqref{eqn:ave_b_3d_2}, and \eqref{eqn:ave_b_3d_3}, this divergence can be shown to vanish:
\begin{equation}
\label{eqn:div_zero_3d}
\nabla \cdot \B^h \Bigl|_{\Tm_{ijk}} = 0.
\end{equation}
Just as in the 2D case, we note that we have actually achieved a globally divergence free
magnetic field, $\B^h(\xv)$.

\section{Numerical examples}
\label{sec:numex}
In this section we apply the proposed 2D and 3D  constrained transport 
schemes to four numerical test cases. For both the 2D and 3D versions,
we consider first a smooth problem to verify the order of accuracy of the proposed
scheme, followed by a problem with shock waves to verify the shock-capturing ability 
of the scheme. All four examples considered in this work
make use of double (2D) or triple (3D) periodic
boundary conditions on the conservative variables:
$\rho$, $\rho \u$, ${\mathcal E}$, and $\B$. 
In our constrained transport scheme, no explicit boundary conditions are needed for $\A$, since 
$\A$ is updated via \eqref{eqn:induction} using velocity and magnetic fields that
satisfy the appropriate boundary conditions.

\subsection{2D smooth Alfv\'en wave problem}
We consider a smooth circular polarized Alfv\'en wave that propagates
in direction $\n = [\cos \phi, \sin \phi, 0]^T$ towards the origin. This problem has been considered
by several authors (e.g., \cite{article:HeRoTa11,article:Ro04b,article:To00}).
The problem consists of  smooth initial data on $\left[0,{1}/{\cos\phi}\right]\times \left[0,{1}/{\sin\phi}\right]$,
where $\phi = \tan^{-1}(0.5)$:
\begin{equation}
\rho = 1, \quad \u = \left[ -u^t \sin \phi, \, \, u^t \cos \phi, \, \, u^r  \right]^T, \quad p = 0.1,
\quad \B = \left[ \cos \phi - u^t \sin \phi, \, \, \sin \phi + u^t \cos \phi, \, \, u^r \right]^T,
\end{equation}
where
\begin{equation}
     u^t = 0.1 \sin \left( 2\pi \left( x \cos \phi +y \sin \phi \right) \right)
     \quad \text{and} \quad
     u^r = 0.1 \cos \left( 2\pi \left( x \cos \phi +y \sin \phi \right) \right).
\end{equation}
This example is used to verify the order of accuracy of the proposed 2D constrained transport scheme.

\begin{figure}[!t]
\begin{center}
\begin{tabular}{cc}
  (a) \includegraphics[width=73mm]{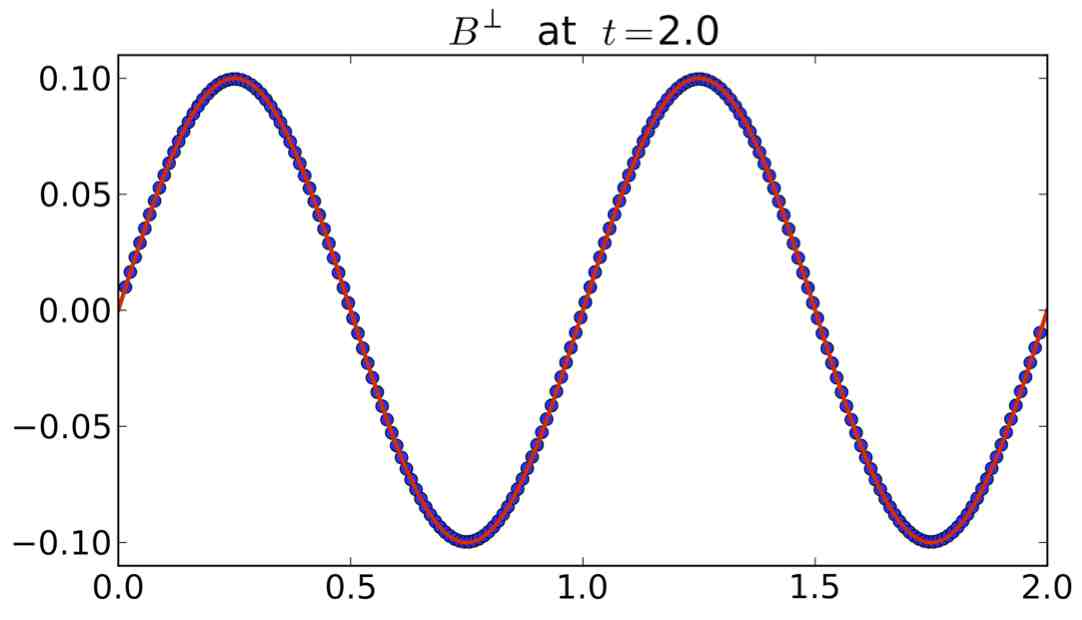} &
  (b) \includegraphics[width=73mm]{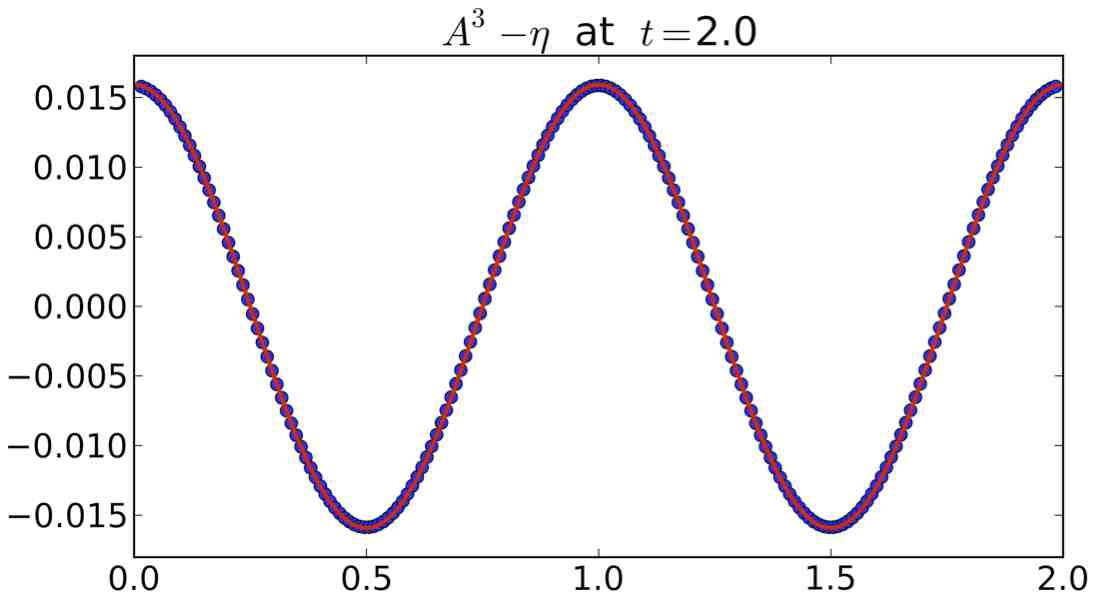} 
\end{tabular}
  \caption{Scatter plots of the constrained transport DG-FEM solution for the 2D smooth Alfv\'en wave on
  	a $16 \times 32$ Cartesian mesh. Panel (a) shows the magnetic field perpendicular to the
	direction of propagation, $B^{\perp} = B^2 \cos \phi - B^1 \sin \phi$, where $\phi = \tan^{-1}(0.5)$,
	as a function of the coordinate along the direction of propagation $\xi = x \cos \phi  + y \sin \phi$.
	Panel (b) shows the periodic part of the magnetic potential, $A^3 - \eta$, where $\eta = y \cos \phi -x \sin \phi$.	
  \label{fig:smooth_alfven}}
\end{center}
\end{figure}

Experimental convergence rates are shown in Table \ref{table:smooth_alfven_2d}. The
errors are calculated by computing the $L^2$-difference between the computed
solution in $\WS^{h}_{2}$ and the exact solution projected into $\WS^{h}_{3}$:
\begin{equation}
	\frac{\| q^{\text{exact}} - q^h \|_{L^2}}{\| q^{\text{exact}} \|_{L^2}} \approx 
		\left[{{\sum_{ij} \left( \sum_{\ell=1}^{6} \left\{ Q^{\text{exact} \, (\ell)}_{ij} - Q^{(\ell)}_{ij}  \right\} 
		+ \sum_{\ell=7}^{10} Q^{\text{exact} \, (\ell)}_{ij} \right)}} \right]^{\half} \cdot
		\left[{{\sum_{ij} \left( \sum_{\ell=1}^{10}  Q^{\text{exact} \, (\ell)}_{ij} \right)}} \right]^{-\half}.
\end{equation}	
In particular, in Table \ref{table:smooth_alfven_2d} we show the computed relative $L^2$-errors 
at time $t=2$ (i.e., after 2 periods of the Alfv\'en wave). Experimental convergence
rates are calculated using a least squares fit through the computed errors. This
table clearly shows the third order convergence of the proposed numerical scheme.
Furthermore, in Figure \ref{fig:smooth_alfven}, we show on a mesh with $16 \times 32$ elements
and at time $t=2$  scatter plots of the computed
(a) magnetic field perpendicular to the direction of propagation ($B^{\perp}$)
and (b) periodic part of the magnetic potential ($A^3 - \left( y \cos\phi -x \sin\phi \right)$). 
These plots show that the computed solution is in very good agreement with the
exact solution.

\subsection{2D Orszag-Tang vortex problem}
Next we consider the Orszag-Tang vortex problem, which is widely considered a standard test example for
MHD in the literature (see T\'oth \cite{article:To00}).
The problem consists of  smooth initial data on $[0,2 \pi] \times [0, 2\pi]$:
\begin{gather}
\rho  = \gamma^2, \quad \u = \Bigl[ -\sin y, \, \sin x, \, 0 \Bigr]^T, \quad
p  = \gamma,  \quad \B = \Bigl[ -\sin y,  \, \sin(2x), \,  0 \Bigr]^T,  \quad A^3 = \frac{1}{2} \cos(2x) + \cos y.
\end{gather}
In this problem, the variable magnetic field eventually causes the smooth initial data to form a strong
rotating shock structure. It has been well documented in the literature (e.g., see T\'oth \cite{article:To00},
 Rossmanith \cite{article:Ro04b}, and Li and Shu \cite{article:LiShu05}) that the formation of this shock structure 
can lead to numerical instabilities in numerical methods that do not control magnetic field
divergence errors. 

\begin{figure}[!t]
\begin{center}
\begin{tabular}{cc}
  (a) \includegraphics[width=70mm]{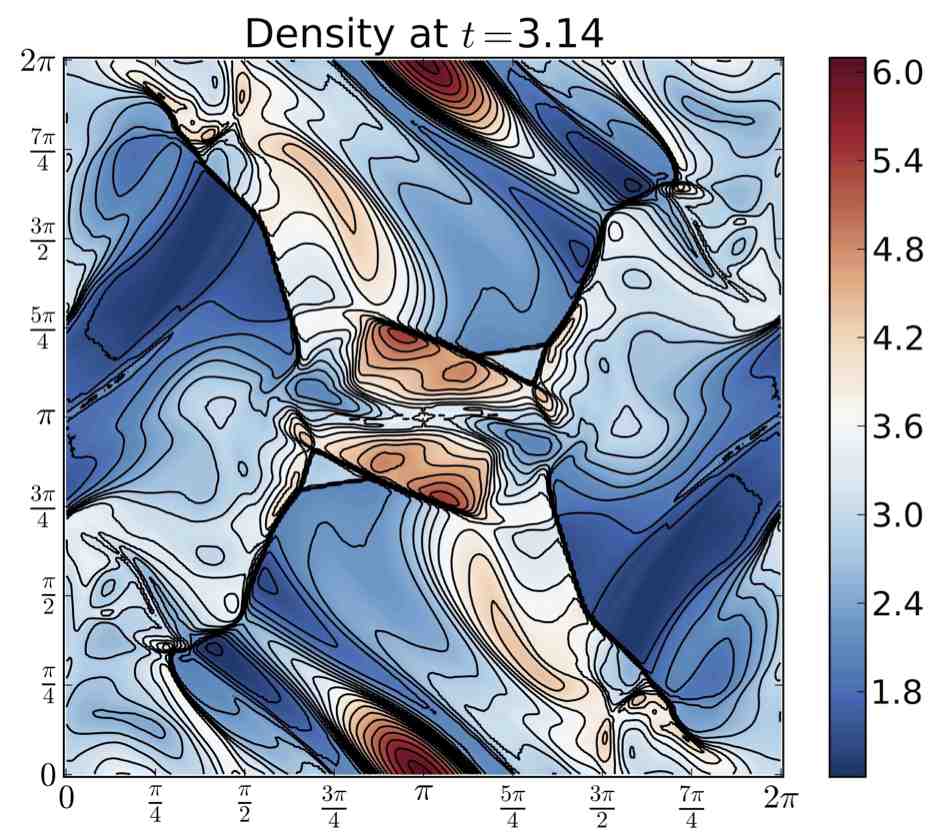} & 
  (b) \includegraphics[width=70mm]{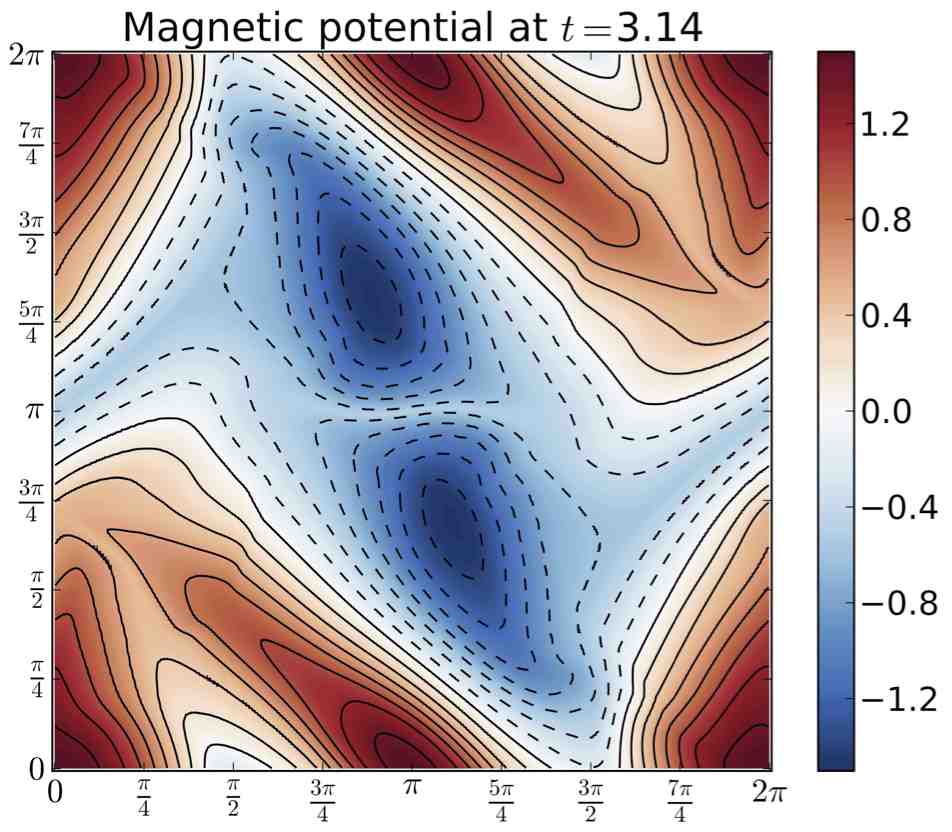}  \\
  (c) \includegraphics[width=70mm]{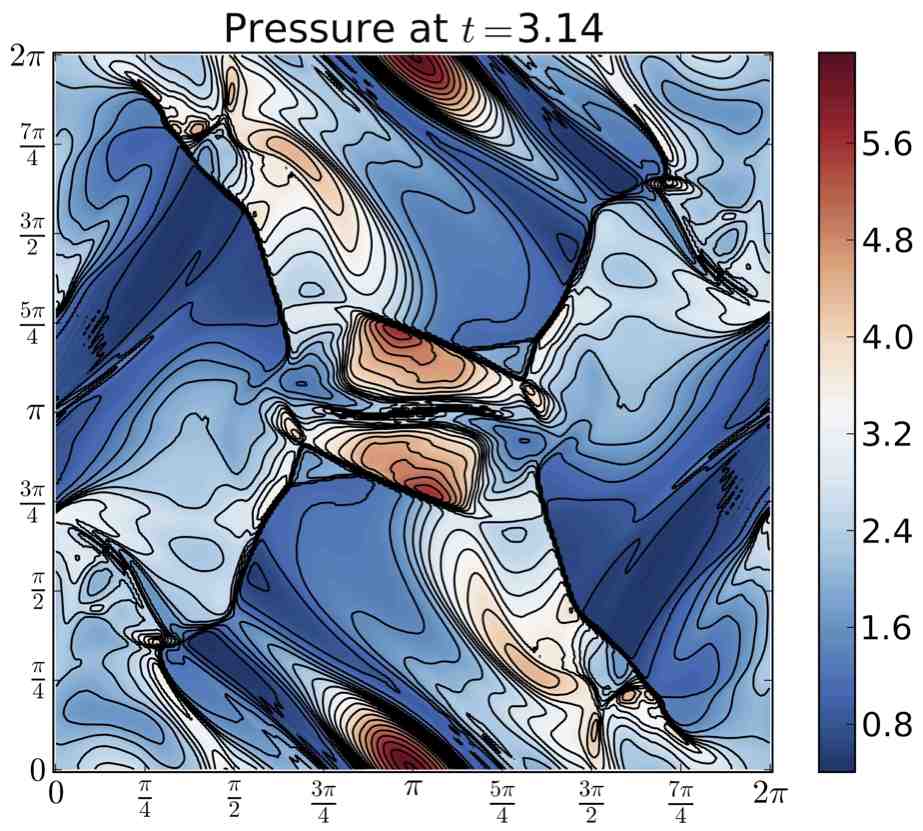} & 
  (d) \includegraphics[width=70mm]{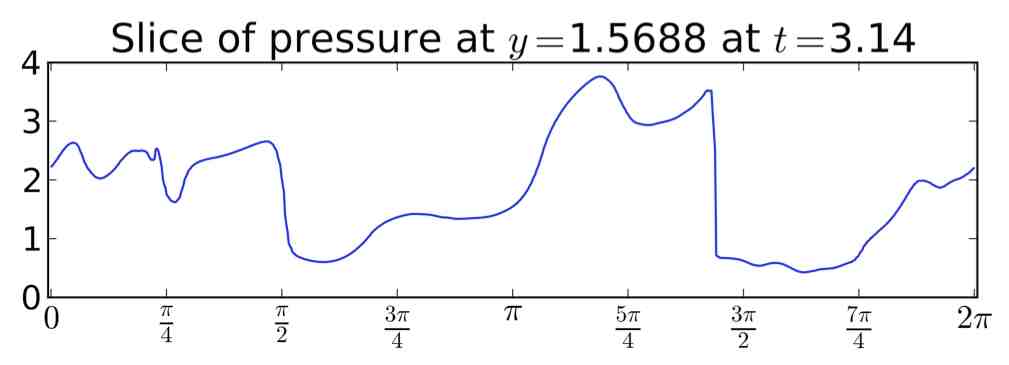}
\end{tabular}
  \caption{Orszag-Tang vortex problem in 2D run on a mesh with $257 \times 257$ mesh elements.
   Shown in these panels are the (a) mass density: $\rho(t,{\bf x})$, (b) magnetic potential: $A^3(t,{\bf x})$,
   (c) pressure: $p(t,{\bf x})$, and (d) a slice of the pressure at $y=1.5688$, all at the final time in
   the solution: $t=3.14$.
  \label{fig:OT-2D}}
\end{center}
\end{figure}

The solution with the proposed scheme is shown in Figure \ref{fig:OT-2D}.
In particular, we show at time $t=3.14$ on a mesh with $257 \times 257$ elements the (a) mass density: $\rho(t,{\bf x})$, (b) magnetic potential: $A^3(t,{\bf x})$, (c) pressure: $p(t,{\bf x})$, and (d) a slice of the pressure at $y=1.5688$.
These results are in good agreement with the published literature, and 
clearly demonstrate the ability of the proposed numerical scheme to remain stable
for a problem with complicated shock structures.

\subsection{3D smooth Alfv\'en wave problem}
In order to verify the order of convergence of the 3D constrained transport
method we consider a 3D version of the smooth Alfv\'en wave problem considered
in Helzel et al. \cite{article:HeRoTa11}.
In this case the wave propagates in the direction
 ${\bf n} = \left[ \cos \phi \cos \theta, \, \sin \phi \cos \theta, \, \sin \theta \right]^T$
towards the origin. Here $\phi$ is an angle with
respect to the $x$-axis in the $xy$-plane and $\theta$ is an angle with respect to the $x$-axis
in the $xz$-plane. We take $\theta = \phi =
\tan^{-1} (0.5) \approx 26.5651^\circ$. 
The problem consists of  smooth initial data on 
$\left[ 0,\left( {\cos \phi \cos \theta} \right)^{-1} \right] \times \left[ 0, \left({\sin \phi \cos
  \theta} \right)^{-1} \right] \times \left[ 0, \left( {\sin \theta} \right)^{-1} \right]$:
\begin{gather}
\rho  = \gamma^2, \quad \u = u^t \, {\bf t} + u^r \, {\bf r}, \quad
p  = \gamma,  \quad \B = \n + B^t \, {\bf t} + B^r \, {\bf r},
\end{gather}
where
\begin{equation}
u^t = 0.1 \sin \left( 2\pi \, {\bf x} \cdot {\bf n} \right), \quad
     u^r = 0.1 \cos \left( 2\pi \, {\bf x} \cdot {\bf n} \right),
\end{equation}
and ${\bf t} = \left[ -\sin \phi, \cos \phi, 0 \right]^T$ and 
${\bf r} = \left[- \cos \phi \sin \theta, -\sin \phi \sin \theta, \cos \theta \right]^T$.
The initial magnetic vector potential is
\begin{align}
 \A &= \left[ 
 	z \sin \phi \cos \theta  - \frac{ \sin \phi
\sin(2 \pi \, {\bf n} \cdot {\bf x})}{20 \pi}, \, \, \,
 x \sin \theta + \frac{ \cos \phi \sin(2 \pi \, {\bf n}
\cdot {\bf x})}{20 \pi}, \, \, \,
  y \cos \phi \cos \theta + \frac{ \cos(2 \pi \, {\bf n}
\cdot {\bf x})}{20 \pi \cos \theta} \right]^T.
\end{align}

Experimental convergence rates are shown in Table \ref{table:smooth_alfven_3d}. The
errors are calculated by computing the $L^2$-difference between the computed
solution in $\WS^{h}_{2}$ and the exact solution projected into $\WS^{h}_{3}$:
\begin{equation}
	\frac{\| q^{\text{exact}} - q^h \|_{L^2}}{\| q^{\text{exact}} \|_{L^2}} \approx 
		\left[{{\sum_{ijk} \left( \sum_{\ell=1}^{10} \left\{ Q^{\text{exact} \, (\ell)}_{ijk} - Q^{(\ell)}_{ijk}  \right\}^2
		+ \sum_{\ell=11}^{20} \left\{ Q^{\text{exact} \, (\ell)}_{ijk} \right\}^2 \right)}} \right]^{\half}
		\cdot \left[{{\sum_{ijk}  \sum_{\ell=1}^{20}  \left\{ Q^{\text{exact} \, (\ell)}_{ijk} \right\}^2 }}\right]^{-\half}.
\end{equation}	
In particular, in Table \ref{table:smooth_alfven_3d} we show the computed relative $L^2$-errors 
at time $t=2$ (i.e., after 2 periods of the Alfv\'en wave). Experimental convergence
rates are calculated using a least squares fit through the computed errors. This
table clearly shows the third order convergence of the proposed numerical scheme 
in all of the magnetic field and magnetic vector potential components.

\subsection{3D Orszag-Tang vortex problem}
Finally, we consider a 3D version of the Orszag-Tang vortex
problem as considered in Helzel et al. \cite{article:HeRoTa11}.
The problem consists of  smooth initial data on $[0,2 \pi] \times [0, 2\pi] \times [0,2\pi]$:
\begin{align}
\rho = \gamma^2, \quad 
\u  = \Bigl[- (1+\varepsilon \sin z)  \sin y, \, \,  (1+ \varepsilon
\sin z) \sin x, \, \, \varepsilon \sin z \Bigr]^T, 
\quad p = \gamma, \quad \B =\Bigl[- \sin y, \, \, \sin(2x), \, \, 0 \Bigr]^T,
\end{align}  
where $\varepsilon = 0.2$.  The initial condition for the magnetic potential is 
\begin{equation}
\A = \Bigl[ 0, \, \, 0,  \, \, \frac{1}{2} \cos(2x) + \cos y \Bigr]^T.
\end{equation}
Just as in the 2D problem, the variable magnetic field eventually causes the smooth initial data to form a strong
rotating shock structure. The formation of this shock structure 
can lead to numerical instabilities in numerical methods that do not control magnetic field
divergence errors. 

\begin{figure}[!t]
\begin{center}
\begin{tabular}{cc}
  (a) \includegraphics[width=70mm]{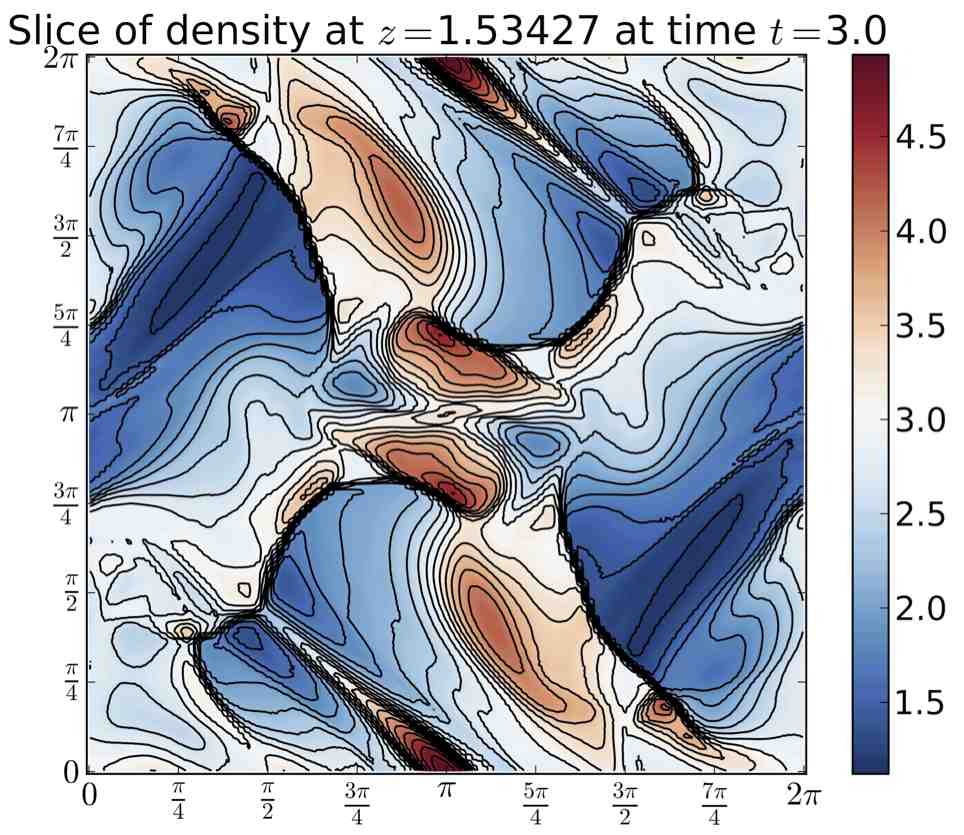} & 
  (b) \includegraphics[width=70mm]{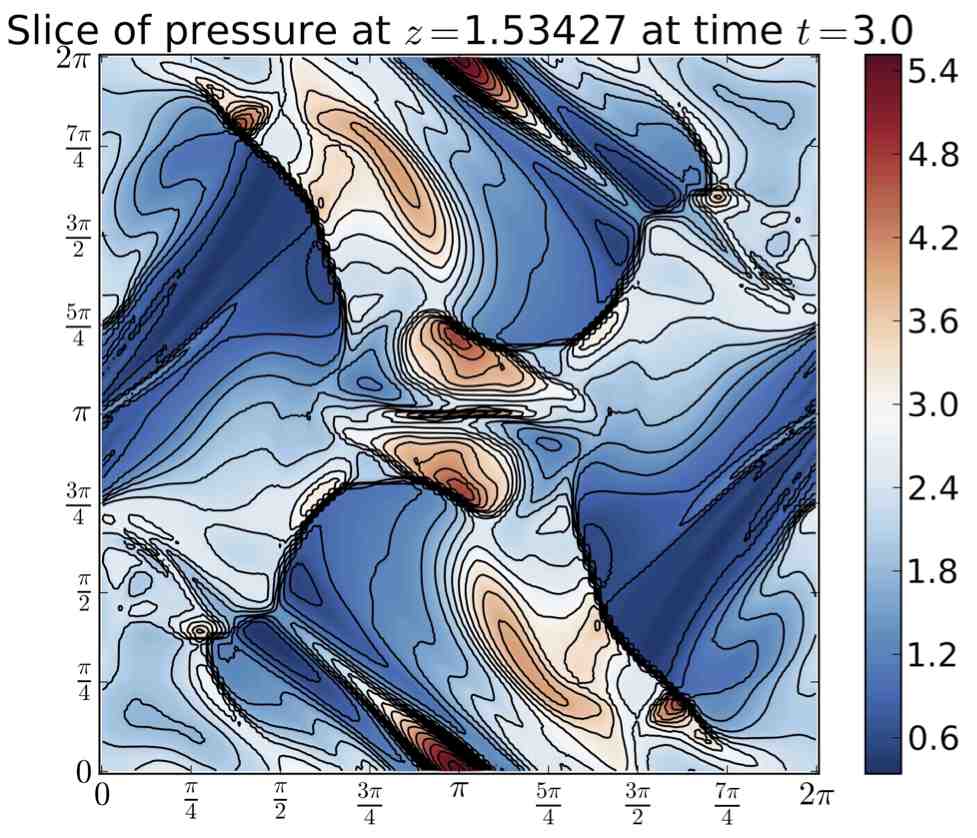} \\
   \multicolumn{2}{c}{(c) \includegraphics[width=70mm]{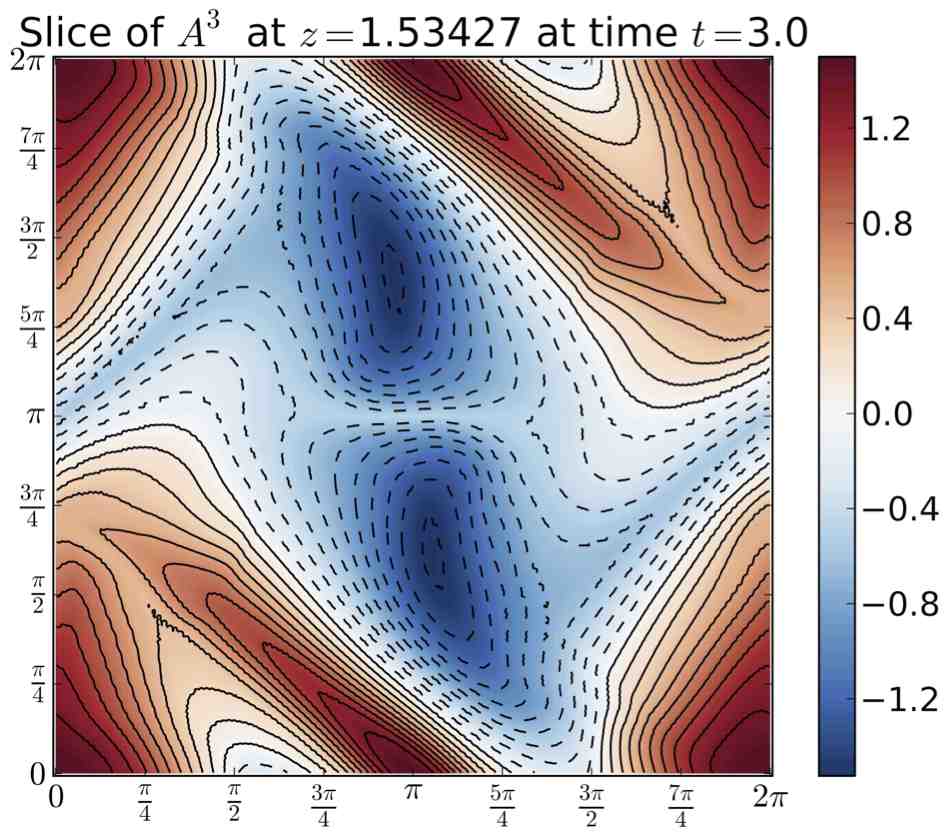}}
\end{tabular}
  \caption{Horizontal slice at $z = z_{32} \approx 1.534266180$ of the 3D Orszag-Tang vortex problem run on a mesh with $129 \times 129 \times 129$ mesh elements.
   Shown in these panels are the (a) mass density: $\rho(t,{\bf x})$, (b) pressure: $p(t,{\bf x})$, and
   (c) $z$-component of the magnetic potential: $A^3(t,{\bf x})$, all at time $t=3$.
  \label{fig:OT-3D-1}}
\end{center}
\end{figure}

\begin{figure}[!t]
\begin{center}
\begin{tabular}{cc}
  (a) \includegraphics[width=70mm]{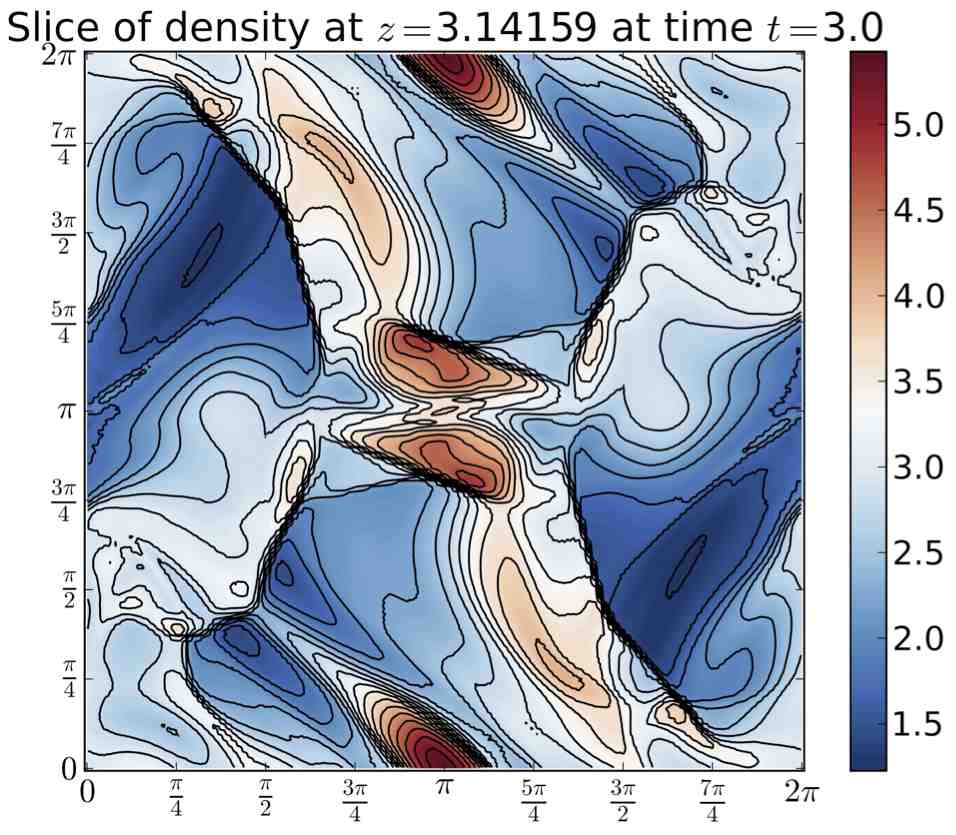} & 
  (b) \includegraphics[width=70mm]{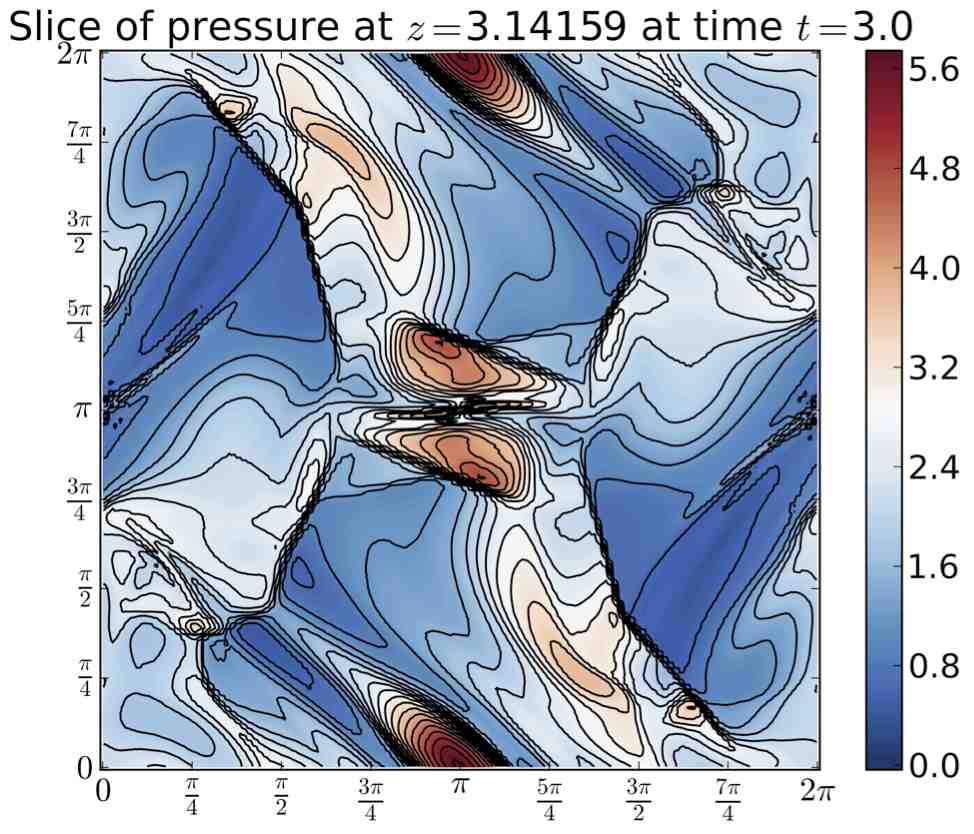} \\
   \multicolumn{2}{c}{(c) \includegraphics[width=70mm]{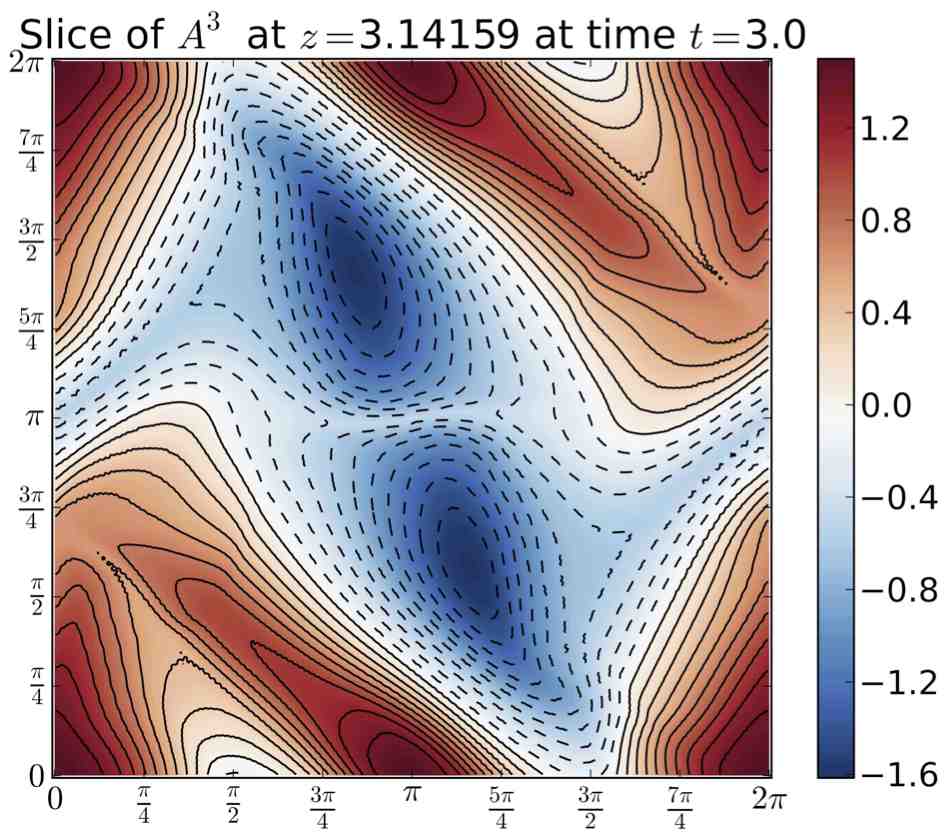}}
\end{tabular}
  \caption{Horizontal slice at $z = z_{65} \approx 3.141592654$ of the 3D Orszag-Tang vortex problem run on a mesh with $129 \times 129 \times 129$ mesh elements.
   Shown in these panels are the (a) mass density: $\rho(t,{\bf x})$, (b) pressure: $p(t,{\bf x})$, and
   (c) $z$-component of the magnetic potential: $A^3(t,{\bf x})$, all at time $t=3$.
  \label{fig:OT-3D-2}}
\end{center}
\end{figure}

\begin{figure}[!t]
\begin{center}
\begin{tabular}{cc}
  (a) \includegraphics[width=70mm]{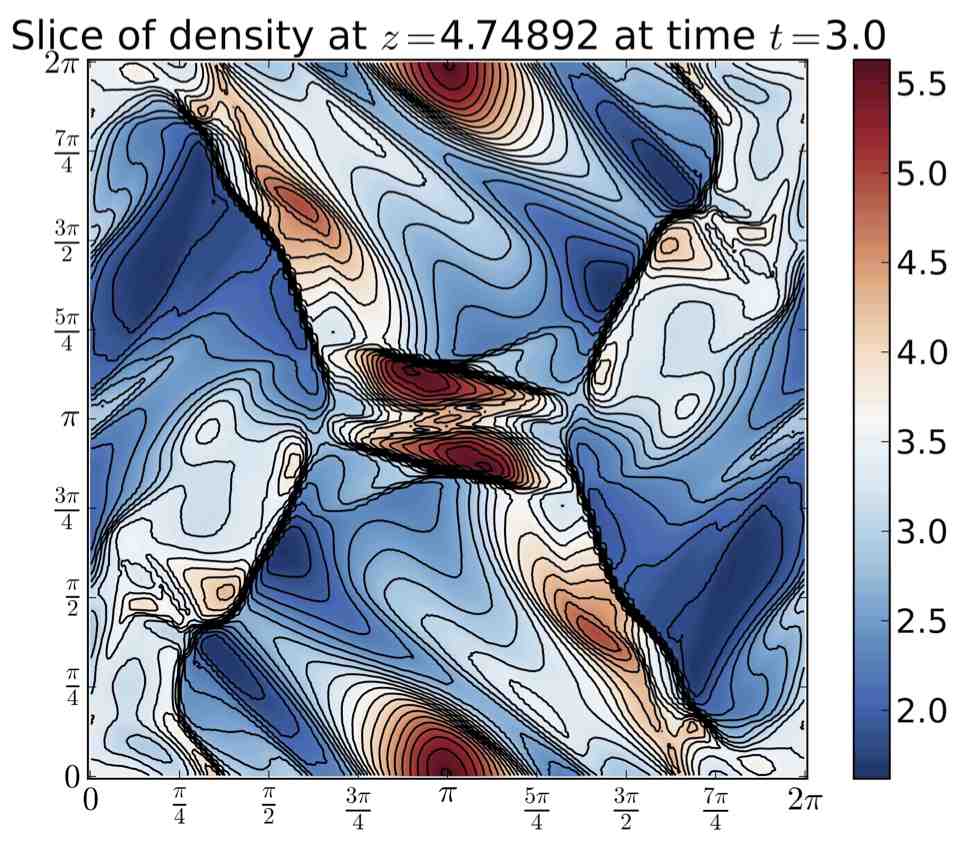} & 
  (b) \includegraphics[width=70mm]{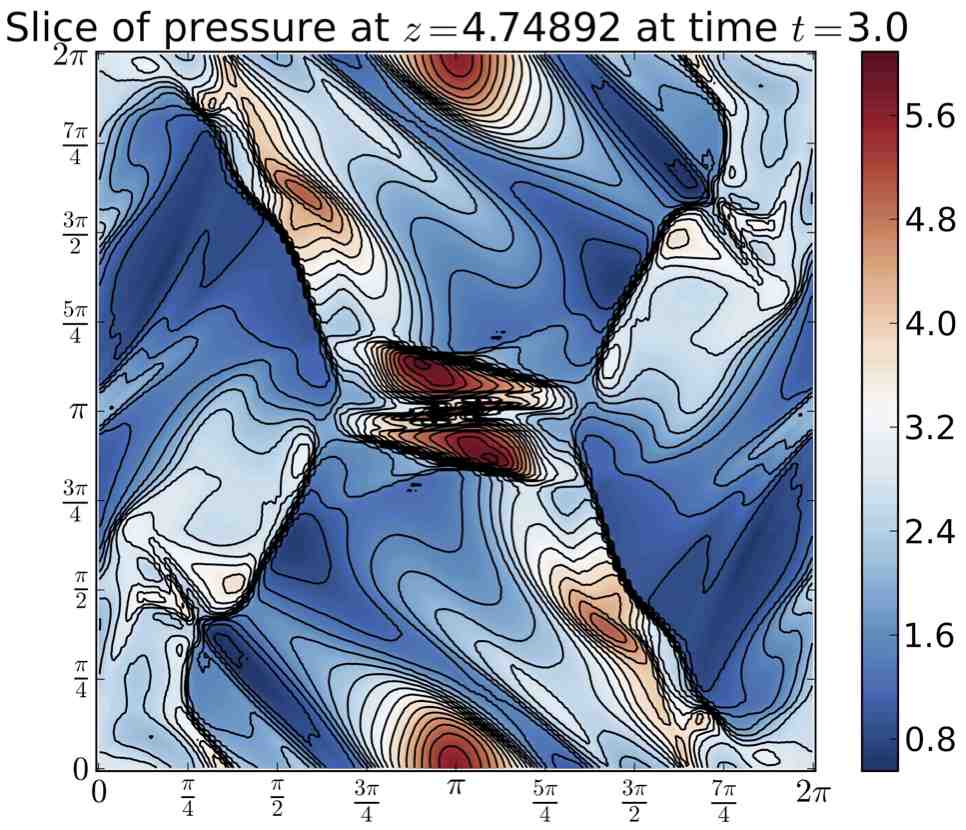} \\
   \multicolumn{2}{c}{(c) \includegraphics[width=70mm]{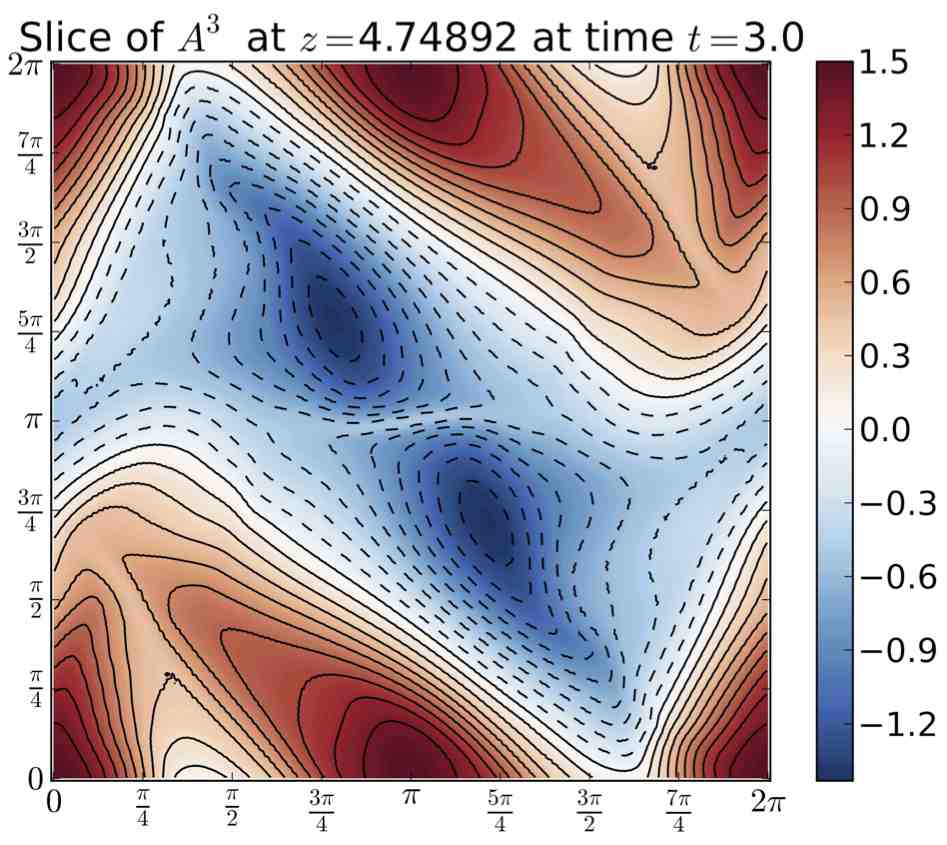}}
\end{tabular}
  \caption{Horizontal slice at $z = z_{98} \approx 4.748919128$ of the 3D Orszag-Tang vortex problem run on a mesh with $129 \times 129 \times 129$ mesh elements.
   Shown in these panels are the (a) mass density: $\rho(t,{\bf x})$, (b) pressure: $p(t,{\bf x})$, and
   (c) $z$-component of the magnetic potential: $A^3(t,{\bf x})$, all at time $t=3$.
  \label{fig:OT-3D-3}}
\end{center}
\end{figure}

We computed a solution on a mesh with $129 \times 129 \times 129$ elements out to time $t=3$.
The results are shown at different horizontal slices 
in Figures \ref{fig:OT-3D-1} ($z  \approx 1.53$),  \ref{fig:OT-3D-2} ($z  \approx 4.75$), and  
\ref{fig:OT-3D-3} ($z \approx 3.14$).
In particular, in each of these plots we show the
(a) mass density: $\rho(t,{\bf x})$, (b) pressure: $p(t,{\bf x})$, and
   (c) $z$-component of the magnetic potential: $A^3(t,{\bf x})$.
These results are in good agreement with those presented in Helzel et al. \cite{article:HeRoTa11}, and 
clearly demonstrate the ability of the proposed numerical scheme to remain stable
for a problem with complicated shock structures.

\begin{table}
\begin{center}
\begin{Large}
\begin{tabular}{|c||c||c||c|}
\hline
{\normalsize {\bf Mesh}} & {\normalsize \bf Rel. $L^2$ error in $B^1$}
& {\normalsize \bf Rel. $L^2$ error in $B^2$}  & {\normalsize \bf Rel. $L^2$ error in $A^3$} \\ 
\hline \hline
{\normalsize $8 \times 16$}   &   {\normalsize $8.8697 \times 10^{-4}$} & {\normalsize $3.4907 \times 10^{-3}$} & {\normalsize $3.3831 \times 10^{-4}$} \\
\hline
{\normalsize $16 \times 32$}   &   {\normalsize $1.0344 \times 10^{-4}$} & {\normalsize $4.0578 \times 10^{-4}$} & {\normalsize $3.8084 \times 10^{-5}$} \\
\hline
{\normalsize $32 \times 64$}   &   {\normalsize $1.3349 \times 10^{-5}$} & {\normalsize $5.2515 \times 10^{-5}$} & {\normalsize $4.6941 \times 10^{-6}$} \\
\hline
{\normalsize $64 \times 128$}   &   {\normalsize $1.7216 \times 10^{-6}$} & {\normalsize $6.798 \times 10^{-6}$} & {\normalsize $5.9009 \times 10^{-7}$} \\
\hline
{\normalsize $128 \times 256$}   &   {\normalsize $2.1876 \times 10^{-7}$} & {\normalsize $8.6605 \times 10^{-7}$} & {\normalsize $7.4004 \times 10^{-8}$} \\
\hline
\hline
{\normalsize {\bf Least squares order}} & {\normalsize 2.9879} & {\normalsize 2.9853} & {\normalsize 3.0329} \\ 
\hline 
\end{tabular} 
\caption{Relative $L^2$ errors in $B^1$, $B^2$, and $A^3$ for several mesh 
resolutions on the 2D smooth Alfv\'en wave problem. For each component we use a 
least squares fit to estimate the order of accuracy. 
\label{table:smooth_alfven_2d}}
\end{Large}
\end{center}
\end{table}

\begin{table}
\begin{center}
\begin{Large}
\begin{tabular}{|c||c||c||c|}
\hline
{\normalsize {\bf Mesh}} & {\normalsize \bf Rel. $L^2$ error in $B^1$}
& {\normalsize \bf Rel. $L^2$ error in $B^2$}  & {\normalsize \bf Rel. $L^2$ error in $B^3$} \\
\hline \hline
{\normalsize $8 \times 16 \times 16$}   &   {\normalsize $1.2159 \times 10^{-3}$} & {\normalsize $6.7737 \times 10^{-3}$} & {\normalsize $3.3457 \times 10^{-3}$} \\
\hline
{\normalsize $16 \times 32 \times 32$}   &   {\normalsize $1.0804 \times 10^{-4}$} & {\normalsize $8.788 \times 10^{-4}$} & {\normalsize $4.9195 \times 10^{-4}$} \\
\hline
{\normalsize $32 \times 64 \times 64$}   &   {\normalsize $1.1507 \times 10^{-5}$} & {\normalsize $1.124 \times 10^{-4}$} & {\normalsize $7.6376 \times 10^{-5}$} \\
\hline
{\normalsize $64 \times 128 \times 128$}   &   {\normalsize $1.5083 \times 10^{-6}$} & {\normalsize $1.4831 \times 10^{-5}$} & {\normalsize $1.0894 \times 10^{-5}$} \\
\hline
\hline
{\normalsize {\bf Least squares order}} & {\normalsize 3.2196} & {\normalsize 2.9472} & {\normalsize 2.7475} \\
\hline 
\end{tabular}

\bigskip

\begin{tabular}{|c||c||c||c|}
\hline
{\normalsize {\bf Mesh}} & {\normalsize \bf Rel. $L^2$ error in $A^1$}
& {\normalsize \bf Rel. $L^2$ error in $A^2$}  & {\normalsize \bf Rel. $L^2$ error in $A^3$} \\
\hline \hline
{\normalsize $8 \times 16 \times 16$}   &   {\normalsize $5.7994 \times 10^{-4}$} & {\normalsize $7.36 \times 10^{-4}$} & {\normalsize $3.5364 \times 10^{-4}$} \\
\hline
{\normalsize $16 \times 32 \times 32$}   &   {\normalsize $8.0689 \times 10^{-5}$} & {\normalsize $9.2977 \times 10^{-5}$} & {\normalsize $4.4725 \times 10^{-5}$} \\
\hline
{\normalsize $32 \times 64 \times 64$}   &   {\normalsize $1.0885 \times 10^{-5}$} & {\normalsize $1.5308 \times 10^{-5}$} & {\normalsize $5.6484 \times 10^{-6}$} \\
\hline
{\normalsize $64 \times 128 \times 128$}   &   {\normalsize $1.4532 \times 10^{-6}$} & {\normalsize $2.2029 \times 10^{-6}$} & {\normalsize $7.2656 \times 10^{-7}$} \\
\hline
\hline
{\normalsize {\bf Least squares order}} & {\normalsize 2.8812} & {\normalsize 2.7755} & {\normalsize 2.9766} \\
\hline
\end{tabular}
\caption{Relative $L^2$ errors in $B^1$, $B^2$, $B^3$, $A^1$, $A^2$, and $A^3$ for several mesh
resolutions on the 3D smooth Alfv\'en wave problem. For each component we use a
least squares fit to estimate the order of accuracy.
\label{table:smooth_alfven_3d}}
\end{Large}
\end{center}
\end{table}

\section{Conclusions}
In this work we showed how to extend the constrained transport framework
originally proposed by Evans and Hawley \cite{article:EvHa88} in the context of the discontinuous Galerkin finite element method (DG-FEM) on both 2D and 3D Cartesian meshes. 
The method presented in this work makes use of two key ingredients:
(1) the introduction of a magnetic vector potential, which is represented
in the same finite element space as the conserved variables, and (2)
the use of a particular divergence-free
reconstruction of the magnetic field, which makes use of the magnetic vector potential
and the predicted magnetic field.
The divergence-free reconstruction presented in this work is slight modification of the reconstruction
method that has been used in other work on ideal MHD. Li, Xu, and Yakovlev \cite{article:LiXuYa11} made
use of this reconstruction in the context of a 2D central DG scheme. Balsara
\cite{article:Ba04} made use this reconstruction in the context of finite volume schemes
and adaptive mesh refinement. 
The novel aspect of this work is that we make direct use of a magnetic vector potential,
thus following in the footsteps of the methods developed in
\cite{article:ChRoTa14,article:HeRoTa11,article:HeRoTa13,article:LoZa04,article:Ro04b,article:De01b}.
An advantage of our approach is that the extension from 2D to 3D is straightforward.
The proposed scheme was then implemented in 2D and 3D using the {\sc DoGPack} software
package and applied on some standard MHD test cases. The results indicate that the
proposed scheme is high-order accurate for smooth problems and is shock-capturing
for problems with complicated multi-dimensional shock structures.

\noindent
{\bf Acknowledgements.}
This work was supported in part by NSF grant DMS--1016202.

\appendix

\section{3D formulas for the construction of a divergence-free magnetic field}
For completeness we include the formulas for computing the 3D divergence-free
reconstruction of the magnetic field in this section. Refer to Figure \ref{fig:yee}(b)
for the positioning of the various magnetic field and magnetic vector potential components.
The magnetic field on left and right faces are defined as follows:
\begin{equation}
\label{eqn:ave_b_3d_1}
\begin{split}
	 b^{1(1)}_{i-\half \, j \, k} &= \frac{1}{\Delta y} \left( {{\mathcal A}^3_{\, i-\half \, j+\half \, k} -   {\mathcal A}^3_{\, i-\half \, j-\half \, k}} \right)
	 	- \frac{1}{\Delta z} \left( {\mathcal A}^2_{\, i-\half \, j \, k+\half} -   {\mathcal A}^2_{\, i-\half \, j \, k-\half} \right),  \\
	 b^{1(2)}_{i-\half \, j \, k} &= \half \left( \Bt^{ \, 1(3)}_{i-1 jk}+\sqrt{3}\Bt^{1(5)}_{i-1 jk}
	 				+\Bt^{ \, 1(3)}_{ijk}-\sqrt{3}\Bt^{ \, 1(5)}_{ijk}\right), \quad
	 b^{1(3)}_{i-\half \, j \, k} = \half \left( \Bt^{ \, 1(4)}_{i-1 jk}+\sqrt{3}\Bt^{1(6)}_{i-1 jk}
	 				+\Bt^{ \, 1(4)}_{ijk}-\sqrt{3}\Bt^{ \, 1(6)}_{ijk}\right), \\
	 b^{1(4)}_{i-\half \, j \, k} &= \half \left( \Bt^{ \, 1(7)}_{i-1 jk} + \Bt^{ \, 1(7)}_{ijk} \right), \quad
	 b^{1(5)}_{i-\half \, j \, k} = \half \left( \Bt^{ \, 1(9)}_{i-1 jk} + \Bt^{ \, 1(9)}_{ijk} \right), \quad
	 b^{1(6)}_{i-\half \, j \, k} = \half \left( \Bt^{ \, 1(10)}_{i-1 jk} + \Bt^{ \, 1(10)}_{ijk} \right).
	 \end{split}
\end{equation}
The magnetic field on front and back faces are defined as follows:
\begin{equation}
\label{eqn:ave_b_3d_2}
\begin{split}
	b^{2(1)}_{i \, j-\half \, k} &= \frac{1}{\Delta z} \left( {\mathcal A}^1_{\, i \, j-\half \, k+\half} - {\mathcal A}^1_{\, i \, j-\half \, k-\half} \right)
	 	- \frac{1}{\Delta x} \left( {\mathcal A}^3_{\, i+\half \, j-\half \, k} -   {\mathcal A}^3_{\, i-\half \, j-\half \, k} \right),  \\
	b^{2(2)}_{i \, j-\half \, k} &=  \half \left( \Bt^{ \, 2(2)}_{i j-1 k}+\sqrt{3}\Bt^{ \, 2(5)}_{i j-1 k}
	 				+\Bt^{ \, 2(2)}_{ijk}-\sqrt{3}\Bt^{ \, 2(5)}_{ijk}\right), \quad
	b^{2(3)}_{i \, j-\half \, k} = \half \left( \Bt^{ \, 2(4)}_{i j-1 k}+\sqrt{3}\Bt^{ \, 2(7)}_{i j-1 k}
	 				+\Bt^{ \, 2(4)}_{ijk}-\sqrt{3}\Bt^{ \, 2(7)}_{ijk}\right), \\
	b^{2(4)}_{i \, j-\half \, k} &= \half \left( \Bt^{ \, 2(6)}_{i j-1 k} + \Bt^{ \, 2(6)}_{ijk} \right), \quad
	b^{2(5)}_{i \, j-\half \, k} = \half \left( \Bt^{ \, 2(8)}_{i j-1 k} + \Bt^{ \, 2(8)}_{ijk} \right), \quad
	b^{2(6)}_{i \, j-\half \, k} = \half \left( \Bt^{ \, 2(10)}_{i j-1 k} + \Bt^{ \, 2(10)}_{ijk} \right).
\end{split}
\end{equation}
The magnetic field on bottom and top faces are defined as follows:
\begin{equation}
\label{eqn:ave_b_3d_3}
\begin{split}
	b^{3(1)}_{i \, j \, k-\half} &= \frac{1}{\Delta x} \left( {\mathcal A}^2_{\, i+\half \, j \, k-\half} - {\mathcal A}^2_{\, i-\half \, j \, k-\half} \right)
	 	- \frac{1}{\Delta y} \left( {\mathcal A}^1_{\, i \, j+\half \, k-\half} -   {\mathcal A}^1_{\, i \, j-\half \, k-\half} \right),  \\
	b^{3(2)}_{i \, j \, k-\half} &= \half \left( \Bt^{ \, 3(2)}_{i j k-1}+\sqrt{3}\Bt^{ \, 3(6)}_{i j k-1}
	 				+\Bt^{ \, 3(2)}_{ijk}-\sqrt{3}\Bt^{ \, 3(6)}_{ijk}\right), \quad
	b^{3(3)}_{i \, j \, k-\half} = \half \left( \Bt^{ \, 3(3)}_{i j k-1}+\sqrt{3}\Bt^{ \, 3(7)}_{i j k-1}
	 				+\Bt^{ \, 3(3)}_{ijk}-\sqrt{3}\Bt^{ \, 3(7)}_{ijk}\right), \\
	b^{3(4)}_{i \, j \, k-\half} &= \half \left( \Bt^{ \, 3(5)}_{i j k-1} + \Bt^{ \, 3(5)}_{ijk} \right), \quad
	b^{3(5)}_{i \, j \, k-\half} = \half \left( \Bt^{ \, 3(8)}_{i j k-1} + \Bt^{ \, 3(8)}_{ijk} \right), \quad
	b^{3(6)}_{i \, j \, k-\half} = \half \left( \Bt^{ \, 3(9)}_{i j k-1} + \Bt^{ \, 3(9)}_{ijk} \right).
\end{split}
\end{equation}

From these normal components on faces, we can reconstruct
a globally divergence-free element-centered definition of the magnetic
field. The 1-component can be written as follows:
\begin{align}
\label{eqn:div-free-b-3d-1}
B^{1(1)}_{ijk} &= \half \left(\aal{1}+\aar{1}\right)+\halfsqt \left( \dxody \left(\bbr{2}-\bbl{2}\right)+ \dxodz \left(\ccr{2}-\ccl{2}\right) \right),
\end{align}
\begin{align}
B^{1(2)}_{ijk} &= \halfsqt \left(\aar{1}-\aal{1}\right)+\halfsqft \left( \dxody \left(\bbr{5}-\bbl{5}\right)
+ \dxodz \left(\ccr{5}-\ccl{5}\right) \right),
\end{align}
\begin{align}
B^{1(3)}_{ijk} &= \half \left(\aar{2}+\aal{2}\right)+\halfsqt \dxodz \left(\ccr{4}-\ccl{4}\right),
\end{align}
\begin{align}
B^{1(4)}_{ijk} &= \half \left(\aar{3}+\aal{3}\right)+\halfsqt \dxody \left(\bbr{4}-\bbl{4}\right),
\end{align}
\begin{align}
B^{1(5)}_{ijk} &= \halfsqt \left(\aar{2}-\aal{2}\right), \quad
B^{1(6)}_{ijk} = \halfsqt \left(\aar{3}-\aal{3}\right),
\end{align}
\begin{align}
B^{1(7)}_{ijk} &= \half \left(\aar{4}+\aal{4}\right), \quad
B^{1(8)}_{ijk} = \halfsqft \left(\dxody \left(\bbl{2}-\bbr{2} \right)+\dxodz \left(\ccl{2}-\ccr{2}\right)\right),
\end{align}
\begin{align}
B^{1(9)}_{ijk} &= \half \left(\aar{5}+\aal{5}\right), \quad
B^{1(10)}_{ijk} = \half \left(\aar{6}+\aal{6}\right),
\end{align}
\begin{align}
B^{1(11)}_{ijk} &= \halfsqft \dxodz \left(\ccl{4}-\ccr{4}\right), \quad 
B^{1(12)}_{ijk} = \halfsqft \dxody \left(\bbl{4}-\bbr{4}\right),
\end{align}
\begin{align}
B^{1(13)}_{ijk} &= \halfsqt \left(\aar{5}-\aal{5}\right), \quad
B^{1(14)}_{ijk} = 0, \quad 
B^{1(15)}_{ijk} = \halfsqt \left(\aar{6}-\aal{6}\right),
\end{align}
\begin{align}
B^{1(16)}_{ijk} &= 0, \quad
B^{1(17)}_{ijk} = \halfsqt \left(\aar{4}-\aal{4}\right),
\end{align}
\begin{align}
B^{1(18)}_{ijk} &= \halfsqtf \left(\dxody \left(\bbl{5}-\bbr{5}\right)+\dxodz \left(\ccl{5}-\ccr{5}\right)\right), \quad
B^{1(19)}_{ijk}= 0, \quad B^{1(20)}_{ijk} = 0.
\end{align}
The 2-component can be written as follows:
\begin{align}
\label{eqn:div-free-b-3d-2}
B^{2(1)}_{ijk} &= \half \left(\bbl{1}+\bbr{1}\right)+\halfsqt \left( 
	\dyodx \left(\aar{2}-\aal{2}\right)+ \dyodz \left(\ccr{3}-\ccl{3} \right) \right), 
\end{align}
\begin{align}
B^{2(2)}_{ijk} &= \half \left(\bbl{2}+\bbr{2} \right),
\end{align}
\begin{align}
B^{2(3)}_{ijk} &= \halfsqt \left(\bbr{1}-\bbl{1}\right)+\halfsqft \left(
	 \dyodx \left(\aar{5}-\aal{5}\right)+ \dyodz \left(\ccr{6}-\ccl{6}\right) \right),
\end{align}
\begin{align}
B^{2(4)}_{ijk} &= \half \left(\bbr{3}+\bbl{3}\right)+\halfsqt \dyodx \left(\aar{4}-\aal{4}\right),
\end{align}
\begin{align}
B^{2(5)}_{ijk} &= \halfsqt \left(\bbr{2}-\bbl{2}\right), \quad
B^{2(6)}_{ijk} = \half \left(\bbr{4}+\bbl{4}\right),
\end{align}
\begin{align}
B^{2(7)}_{ijk} &= \halfsqt \left(\bbr{3}-\bbl{3}\right), \quad
B^{2(8)}_{ijk} = \half \left(\bbr{5}+\bbl{5}\right),
\end{align}
\begin{align}
B^{2(9)}_{ijk} &= \halfsqft \left(\dyodx \left(\aal{2}-\aar{2}\right)+\dyodz \left(\ccl{3}-\ccr{3}\right)\right),
\end{align}
\begin{align}
B^{2(10)}_{ijk} &= \half \left(\bbl{6}+\bbr{6}\right), \quad
B^{2(11)}_{ijk} = \halfsqt \left(\bbr{5}-\bbl{5}\right), \quad
B^{2(12)}_{ijk} = B^{2(13)}_{ijk} = 0,
\end{align}
\begin{align}
B^{2(14)}_{ijk} &= \halfsqft \dyodx \left(\aal{4}-\aar{4}\right), \quad
B^{2(15)}_{ijk} = 0, \quad 
B^{2(16)}_{ijk} = \halfsqt \left(\bbr{6}-\bbl{6}\right),
\end{align}
\begin{align}
B^{2(17)}_{ijk} &= \halfsqt \left(\bbr{4}-\bbl{4}\right), \quad 
B^{2(18)}_{ijk} = 0,
\end{align}
\begin{align} 
B^{2(19)}_{ijk} &= \halfsqtf \left(\dyodx \left(\aal{5}-\aar{5}\right)+\dyodz \left(\ccl{6}-\ccr{6}\right)\right), 
\quad B^{2(20)}_{ijk} = 0.
\end{align}
The 3-component can be written as follows:
\begin{align}
\label{eqn:div-free-b-3d-3}
B^{3(1)}_{ijk} &= \half \left(\ccl{1}+\ccr{1}\right)+\halfsqt \left(
  \dzodx \left(\aar{3}-\aal{3}\right)+ \dzody \left(\bbr{3}-\bbl{3} \right) \right),
\end{align}
\begin{align}
B^{3(2)}_{ijk} &= \half \left(\ccr{2}+\ccl{2}\right), \quad 
B^{3(3)}_{ijk} = \half \left(\ccr{3}+\ccl{3}\right),
\end{align}
\begin{align}
B^{3(4)}_{ijk} &= \halfsqt \left(\ccr{1}-\ccl{1}\right)+\halfsqft \left(\dzodx \left(\aar{6}-\aal{6}\right)+\dzody \left(\bbr{6}-\bbl{6} \right)\right),
\end{align}
\begin{align}
B^{3(5)}_{ijk} &= \half \left(\ccr{4}+\ccl{4}\right), \quad
B^{3(6)}_{ijk} = \halfsqt \left(\ccr{2}-\ccl{2}\right),
\end{align}
\begin{align}
B^{3(7)}_{ijk} &= \halfsqt \left(\ccr{3}-\ccl{3}\right),  \quad
B^{3(8)}_{ijk} = \half \left(\ccr{5}+\ccl{5}\right),
\end{align}
\begin{align}
B^{3(9)}_{ijk} &= \half \left(\ccr{6}+\ccl{6}\right),
\end{align}
\begin{align}
B^{3(10)}_{ijk} &= \halfsqft \left(\dzodx \left(\aal{3}-\aar{3}\right)+\dzody \left(\bbl{3}-\bbr{3}\right)\right),
\end{align}
\begin{align}
B^{3(11)}_{ijk} &= 0, \quad
B^{3(12)}_{ijk} = \halfsqt \left(\ccr{5}-\ccl{5}\right), \quad
B^{3(13)}_{ijk} = 0,
\end{align}
\begin{align}
B^{3(14)}_{ijk} &= \halfsqt \left(\ccr{6}-\ccl{6}\right),
\end{align}
\begin{align}
B^{3(15)}_{ijk} &= B^{3(16)}_{ijk} = 0, \quad 
B^{3(17)}_{ijk} = \halfsqt \left(\ccr{4}-\ccl{4}\right), \quad 
B^{3(18)}_{ijk} = B^{3(19)}_{ijk} = 0,
\end{align}
\begin{align}
\label{eqn:div-free-b-3d-3-end}
B^{3(20)}_{ijk} &= \halfsqtf \left(\dzodx \left(\aal{6}-\aar{6}\right)+\dzody \left(\bbl{6}-\bbr{6}\right)\right).
\end{align}

\end{document}